\documentclass[a4,11pt]{article}

\usepackage[latin1]{inputenc}
\usepackage{color}
\usepackage{latexsym}
\usepackage{amssymb}
\usepackage{graphicx}
\usepackage{enumerate}
\usepackage{extarrows}
\newtheorem{Theorem}{Theorem}[part]

\newtheorem{Proposition}{Proposition}[part]
\newtheorem{Assumption}{Assumption}[part]
\newtheorem{Lemma}{Lemma}[part]

\newtheorem{Remark}{Remark}[part]

\def \s{~~~}
\def \nn {\nonumber}

\def \Int{\displaystyle\int}

\def \R{\mathbb{R}}

\def \E{\mathbb{E}}

\def \eps{\varepsilon}

\def \ep{\hbox{ }\hfill$\Box$}
\def \2 {\vspace{2mm}}

\def\beqs{\begin{eqnarray*}}
\def\enqs{\end{eqnarray*}}
\def\beq{\begin{eqnarray}}
\def\enq{\end{eqnarray}}

\begin{document}

\title{Discrete time approximation of decoupled Forward-Backward SDE  driven by  pure jump Lévy-processes}

\author{Soufiane Aazizi \thanks{ Department of Mathematics, Faculty of Sciences Semlalia
Cadi Ayyad University, B.P. 2390 Marrakesh, Morocco. \sf aazizi.soufiane@gmail.com}}


\maketitle

\begin{abstract}
We present a new algorithms to discretize a decoupled forward backward stochastic differential equations
driven by  pure jump Lévy process (FBSDEL in short). The method is built  in two steps. Firstly, we approximate  the FBSDEL by a forward backward stochastic differential equations driven by a Brownian motion and Poisson process (FBSDEBP in short), in which we replace the small jumps by a Brownian motion. Then, we prove the convergence of the approximation when the size of small jumps $\eps$ goes to $0$. In the second step, we obtain the $L^p$ Hölder continuity of the solution of FBSDEBP and we construct two numerical schemes for this FBSDEBP. Based on the $L^p$ Hölder estimate, we prove the convergence of the scheme when the number of time steps $n$ goes to infinity.  Combining these two steps leads to prove the convergence of numerical schemes to the solution of FBSDEL.
\end{abstract}


\vspace{13mm}

\noindent {\bf Key words~:} Discrete-time approximation,  Euler scheme, decoupled forward-backward SDE with jumps, Small jumps,  Malliavin calculus.

\vspace{5mm}

\noindent {\bf MSC Classification:} 60H35, 60H07, 60J75

\newpage

\section {Introduction and summary}
In this paper, we are concerned by discretization of a system of decoupled  forward-backward stochastic differential equation (FBSDEs in short) driven by a pure jump Lévy process
\begin{equation}
\label{FBSDE_L}
\left\{
\begin{array}{lll}
X_t&=&X_0 + \int_0^t b(X_r)dr + \int_0^t \int_\mathbb{R}\beta(X_{r^-})\bar{M}(de,dr),  \\
Y_t&=& g(X_T) +\int_t^Tf(\Theta_r)dr - \int_t^T \int_\mathbb{R}V_r \bar{M}(de,dr)
\end{array}
\right.
\end{equation}
Here $\Theta :=\left(X, Y, \int_E \rho(e)Ve \nu(de)\right)$ and  $\bar{M}(E, t  )  = \int _{E \times [0, t] } e \bar{\mu}(de,dr)$  where
 $\bar{\mu}(de,dr) :=\mu(de,dr) - \nu(de)dr$ an independent compensated Poisson measure and $ \mu$  a Poisson random measure on $\R\times [0, T]$ with intensity $\nu$ satisfying $\int 1 \wedge |e|^2\nu(de) < \infty $.
\2 \\
Numerical discretization schemes for FBSDE have been studied by many authors. In the no-jump case,  Ma et al. \cite {MPY94} developed the first step algorithm to solve a class of general forward-backward SDE. Douglas et al. \cite{DMP96} suggest  a finite difference approximation of the associated PDE. Other discrete scheme have been considered in \cite{BDM01}, \cite{C97} and \cite{CMM98} mainly based on approximation of the Brownian motion by some discrete process, Gobet et al. \cite{GLW06} proposed an adapted Longstaff and Schwartz algorithm  based on non-parametric regressions. In the jump case, to our knowledge, there is only the work of Bouchard and Elie \cite{BE08} in which the authors propose a Monte-Carlo methods in the case when $\nu(\R) < \infty$.
\2 \\
The main motivation to study the numerical scheme of a systems of above form, is to treat the case when $\nu(\R)=\infty$, which means the existence of an infinite  number of jumps in every interval of non-zero length a.s.. In this sense, we should mention the important work on the approximation of stochastic differential equation studied by Kohatsu-Higa and Tankov \cite{KT10}.
\2 \\
Since we are interested in the case of $\nu(\R) =\infty$, we will follow the idea  of \cite{KT10} to approximate (\ref{FBSDE_R}) without cutoff the small jumps smaller than $\eps$, which should improve the approximation scheme. Then by using the approximation result of Asmussen and Rosinski \cite{AR01} we replace the small jumps of the driven-Lévy process  with  $\sigma(\eps) W$ where $W$ is a standard Brownian motion and
$ \label{sigma}\sigma^2(\eps):=\int_{E^\eps} e^2 \nu(de)$.
\2 \\
In the aim to approximate (\ref{FBSDE_L}), we cut the jumps at $\eps$ as the following
\begin{equation}
\label{FBSDE_R}
\left\{
\begin{array}{lll}
X_t &=& X_0+ \int_0^t b(X_r)dr + \int_0^t \beta(X_{r^-})dR_r + \int_0^t \int_{E_\eps }  \beta(X_{r^-})\bar{M}(de,dr)\\
Y_t &=& g(X_T) + \int_t^Tf(\Theta_r)dr -\int_t^T  V_rdR_r -\int_t^T \int_{E_\eps }  V_r \bar{M}(de,dr)
\end{array}
\right.
\end{equation}
where $R_t = \int_0^t \int_{|e|\le \eps}  e \bar{M}(de,dr)$, $E^\eps :=\{e\in \R, \mbox{ s.t } / \ |e|\leq \eps\}$, $E_\eps :=\{e\in \R, \mbox{ s.t } / \ |e| > \eps\}$ and  $E:=\R= E^\eps  \cup E_\eps$.
\2 \\
The idea we propose is to  discretize the solution of (\ref{FBSDE_L}) in two steps. In the first step, we approximate (\ref{FBSDE_R}) by the following FBSDE:
\begin{equation}
\label{FBSDE_Eps}
\left\{
\begin{array}{lll}
X_t^\eps&=& X_0^\eps+ \int_0^t b(X_r^\eps)dr + \int_0^t \beta(X_r^\eps)\sigma(\eps)dW_r + \int_0^t \int_{E_\eps}\beta(X_{r^-}^\eps) \bar{M}(de,dr)\\
Y_t^\eps &=& g(X_T^\eps) + \int_t^Tf(\Theta_r^\eps)dr -\int_t^T  Z_r^\eps dW_s -\int_t^T \int_{E_\eps} U_r^\eps  (e) \bar{M}(de,dr)
\end{array}
\right.
\end{equation}
Here $\Theta^\eps :=\Big(X^\eps, Y^\eps, \Gamma^\eps \Big)$ and $ \Gamma^\eps := \int_{E_\eps}\rho(e)U^\eps(e)e \nu(de)$. Further, we  show that for a finite measure $m$ defined by $m(E) := \int_E e^2 \nu(de)$, our error
\beqs Err^2_\eps(Y, V)&:=& \mathbb{E}\left[\sup_{t\leq T}|Y_t  - Y_t^\eps|^2 \right] + \mathbb{E} \left[\sup_{t\leq T} \left| \int_0^tV_r dR_r - \int_0^tZ^\eps_rdW_r \right|^2 \right] \\
    && + \mathbb{E}\left[\int_0^T\int_{E_\eps} |V_r - U_r^\eps (e)|^2m(de)dr \right],
\enqs
is controlled by $\sigma(\eps)^2$, which means that the solution of (\ref{FBSDE_Eps}) converges to the solution of (\ref{FBSDE_L}), as the size of small jumps $\eps$ goes to $0$ (See Remark \ref{RemYeps}). We also derive the upper bound
\beq \mathbb{E}\left[\sup_{t\leq T}|X_t  - X_t^\eps|^2 \right] \leq C \sigma(\eps)^2 . \enq
The second step consists of discretizating  the approximated FBSDE (\ref{FBSDE_Eps}) and studying its convergence to (\ref{FBSDE_R}). For this purpose we consider two numerical schemes, the first one is based on discrete-time approximation of decoupled FBSDE derived by  Bouchard and Elie \cite{BE08}. More precisely, for a fixed $\eps$, given a regular grid
$\pi = \{t_i=iT/n, i=0, 1, ..., n.\}$, the authors  approximate $X^\eps$ by its Euler scheme $\bar{X}^\pi$ and $(Y^\eps, Z^\eps, \Gamma^\eps)$ by the
discrete-time process $(\bar{Y}^\pi_{t},\bar{Z}^\pi_{t},\bar{\Gamma}^\pi_{t} )$
\begin{equation}
\label{Approximation}
\left\{
\begin{array}{lll}
\bar{X}^\pi_{t_i+1} &=& \bar{X}^\pi_{t_i}+\frac{1}{n}b(\bar{X}^\pi_{t_i}) + \beta(\bar{X}^\pi_{t_i})\sigma(\eps)\Delta W_{i+1} + \int_{E_\eps}\beta(\bar{X}^\pi_{t_i}) \bar{M}(de,(t_i, t_{i+1}] ) \\
\bar{Z}^\pi_{t}&=&n\mathbb{E} \left[ \bar{Y}^\pi_{t_{i+1}} \Delta W_{i+1} / \mathcal{F}_{t_i} \right] \\
\bar{\Gamma}^\pi_{t}&=&n\mathbb{E} \left[ \bar{Y}^\pi_{t_{i+1}} \int_{E_\eps} \rho(e) \bar{M} (de,(t_i, t_{i+1}] )   / \mathcal{F}_{t_i} \right]\\
\bar{Y}^\pi_{t}&=&\mathbb{E} \left[ \bar{Y}^\pi_{t_{i+1}}  / \mathcal{F}_{t_i} \right] + \frac{1}{n} f(\bar{X}^\pi_{t_i}, \bar{Y}^\pi_{t_i}, \bar{\Gamma}^\pi_{t_i} )
\end{array}
\right.
\end{equation}
on each interval $[t_i, t_{i+1})$, where the terminal value $\bar{Y}^\pi_{t_n}:=g(\bar{X}^\pi_{t_n}) $. Under Lipschitz continuity of the solution, the authors  proved that the
discretization error
 \begin{eqnarray} \label{errorYPhi}
\overline{Err}^2_n(Y^\eps, Z^\eps, \Gamma^\eps) &:=&  \sup_{t\leq T} \mathbb{E} \left[ |Y^\eps_t - \bar{Y}^\pi_t|^2 \right] + \int_0^T\E\big[|Z^\eps_t - \bar{Z}^\pi_t|^2 +  |\Gamma^\eps_t - \bar{\Gamma}^\pi_t |^2\big]dt
\end{eqnarray}
achieves the optimal convergence rate $n^{-1/2} $.
Finally, we derive the first main result of this paper in Proposition \ref{PropErrNEps} showing that the approximation-discretization error
\begin{eqnarray}
\overline{Err}^2_{(n,\eps)}(Y, V) &:=& \sup_{t\leq T} \mathbb{E} \left[ |Y_t - \bar{Y}^\pi_t|^2 + \Big|\int_0^t V_rdR_r - \int_0^t\bar{Z}_r^\pi dW_r \Big|^2 \right]+ \|\Gamma - \bar{\Gamma}^\pi \|^2 _{H^2},
\end{eqnarray}
is bounded by $C (n^{-1} + \sigma(\eps)^2 )$ and converges to $0$ as $(n, \eps)$ tends to $(\infty,0)$, where  $\Gamma := \int_{E_\eps} \rho(e)Ve \nu(de)$ . Taking $\eps = n^{-1/2}$, our approximation-discretization achieves the optimal convergence rate $n^{-1/2}$.
\vspace{2mm}\\
The second numerical scheme has been inspired from the paper of Hu, Nualart and Song \cite{HNS11}. Where the authors study a backward stochastic differential equation driven by a Brownian motion with general terminal variable $\xi$. They propose a new scheme using the  representation of $Z^\eps$ as the trace of the Malliavin derivatives of $Y^\eps$. Their discretization scheme is based on the $L^p$-Hölder continuity of the solution $Z^\eps$, to obtain an estimate of the form
\beqs
\mathbb{E} |Z^\eps_t - Z^\eps_s|^p \leq K|t-s|^\frac{p}{2},
\enqs
which implies the existence of a $\gamma$-Hölder continuous version of the process $Z^\eps$ for any $\gamma < \frac{1}{2} - \frac{1}{p}$.
In this sense, our article extend the work done in  \cite{HNS11} to a forward-backward stochastic differential equation with jumps and terminal
value $g(X^\eps_T)$. Similarly to \cite{HNS11}, we obtain the following regularity of $\Gamma^\eps$
\beqs \mathbb{E}|\Gamma^\eps_t - \Gamma^\eps_s|^p &\leq& C |t-s|^\frac{^p}{2}, \label{Gt-Gs}\enqs
which allows us to deduce the existence of a $\gamma$-Hölder continuous version of the process $\Gamma^\eps$ for any $\gamma < \frac{1}{2} - \frac{1}{p}$.
Finally, on one hand, we use the representation of $Z^\eps$ and $\Gamma^\eps$ as the trace of Malliavin derivative of $Y^\eps$ to derive our a new extended discretization scheme for the solution $(Y^\eps, Z^\eps, \Gamma^\eps)$ of (\ref{FBSDE_Eps}). From other hand we approximate $X^\eps$ by $X^\pi $ the continuous-time version of the Euler  scheme, that is for a fixed $\eps >0$ \\
\begin{equation}\label{ApprxFBSDE}\left\{
\begin{array}{lll}
X^\pi_t&=& X^\pi_{\phi^n_t} + \sigma(\eps)b(X^\pi_{\phi^n_t})( t -\phi^n_t)+  \sigma(\eps)\beta(X^\pi_{\phi^n_t})(W_t - W_{\phi^n_t}) + \int_{E_\eps} \beta(X^\pi_{\phi^n_t})\bar{M}(de,(t, \phi^n_t]).\\
Y^\pi_{t_i}&=&\mathbb{E}\Big[ Y^\pi_{t_{i+1}} + f(\Theta^\pi_{t_{i+1}})\Delta t_{i+1}/ \mathcal{F}_{t_i}\Big]\\
Z^\pi_{t_i}&=&\mathbb{E}\Big[\mathcal{E}^\pi_{t_{i+1}, t_n}\partial_x g(X^\pi_T) D_{t_i}X^\pi_T +\sum_{k=i}^{n-1}\mathcal{E}^\pi_{t_{i+1}, t_{k+1}} \partial_xf(\Theta^\pi_{t_{k+1}})D_{t_i}X^\pi_{t_k}\Delta t_{k}/ \mathcal{F}_{t_i}\Big]\\
\Gamma^\pi_{t_i}&=&\E\Big[\int_{E_\eps}\rho(e)\big[\mathcal{E}^{e,\pi}_{t_{i+1}, t_n}D_{t_i,e}g(X^\pi_T) +\sum_{k=i}^{n-1}\mathcal{E}^{e,\pi}_{t_{i+1}, t_{k+1}} \alpha^\pi_{t_i, t_{k+1}}D_{t_i,e}X^\pi_{t_{k+1}}\Delta t_k\big]\nu(de)\big/ \mathcal{F}_{t_i}\Big]\\
\end{array}\right.
\end{equation}
with terminal values $Y^\pi_{t_n}=g(X^\pi_T)$,   $ Z^\pi_{t_n}=\sigma(\eps) \partial_x g(X^\pi_T)\beta(X^\pi_T) $ and  $ U^\pi_{t_n,e}=g(X^\pi_T + \beta(X^\pi_T)) - g(X^\pi_T)$,
where $\phi^n_t$,  $ \mathcal{E}^\pi_{t_i, t_j}$ and $ \mathcal{E}^{e,\pi}_{t_i, t_j}$ are detailed in section 4. \2 \\
The key-ingredient for computation of discretization error, is based on the $L^p$-Hölder continuity of the solution $(Y^\eps,Z^\eps,\Gamma^\eps)$.
This allows us to prove that \\
\beqs
Err^2_n(Y^\eps,Z^\eps,\Gamma^\eps)&:=& \mathbb{E} \max_{0 \leq i \leq n} \Big[|Y^\eps_{t_i} -Y^\pi_{t_i}|^2 + |Z^\eps_{t_i} -Z^\pi_{t_i}|^2 + |\Gamma^\eps_{t_i} -\Gamma^\pi_{t_i}|^2  \Big],
\enqs
is controlled by  $|\pi|^{1-\frac{1}{log\frac{1}{|\pi|}}} $. Then we obtain the second main result of this article in Theorem \ref{ThErrorYEpsi}, which stating that
\beqs
Err^2_{n, \eps}(Y, V) &:=& \max_{0 \leq i \leq n} \sup_{t\in[t_i, t_{i+1}]} \mathbb{E}  \left[|Y_t -Y^\pi_{t_i}|^2 \right] + \E\left|\int_0^TV_rdR_r  -\sum_{i=0}^{n-1}Z^\pi_{t_i}\Delta W_{t_i}\right|^2  \\
                                && + \sum_{i=0}^{n-1}\int_{t_i}^{t_{i+1}}\mathbb{E}|\Gamma_t -\Gamma^\pi_{t_i}|^2 dt,
\enqs
is of the order $\sigma(\eps)^2  +  |\pi|^{1-\frac{1}{log\frac{1}{|\pi|}}}$ and converges to 0 as the discretization step $(\eps, n) $ tends to $(0,\infty)$.
\vspace{2mm}\\
The importance of the above scheme, is it can be adapted to the case when the a terminal value is not given by the forward diffusion equation $X^\eps$, as it 's the case in \cite{HNS11}. However, this scheme remains to be further investigated.
\2 \\
The two numerical schemes above are not directly implemented in practice and require an important procedure to simulate the conditional expectation.
However, there exist different technics which can be adapted to our setting to compute this conditional expectation and we shall only mention the papers:
\cite{BP03}, \cite{BN04} , \cite{CLP02} and \cite{GLW06}. \2 \\
The paper is organized as follows. In Section 2, we prove the convergence of the approximated scheme. In Section 3, we describe discrete-time scheme introduced in \cite{BE08} and state our first main convergence result. In section 4, we extend the new discrete scheme of \cite{HNS11} and state our second main result. We also discuss a general case of BSDE. Section 5, is devoted to Malliavin calculus for a class of FBSDE with jumps, we then get the
  $L^p$-Hölder continuity of $Z^\eps$ and $\Gamma^\eps$ via the trace of the Malliavin derivatives of $Y^\eps$.
\section{Approximation of decoupled FBSDE driven by pure jump Lévy processes }
\setcounter{equation}{0} \setcounter{Assumption}{0}
\setcounter{Theorem}{0} \setcounter{Proposition}{0}
\setcounter{Corollary}{0} \setcounter{Lemma}{0}
\setcounter{Definition}{0} \setcounter{Remark}{0}

Let $(\Omega, \mathcal{F}, \mathbb{F}=(\mathcal{F})_{t \leq T}, \mathbb{P})$ be a stochastic basis such that
$\mathcal{F}_0$ contains the $\mathbb{P}$-null sets, $\mathcal{F}_T=\mathcal{F}$ and $\mathbb{F}$ satisfies the usual assumptions.
We assume that $\mathbb{F}$ is generated by a one-dimensional Brownian motion $W$ and an independent Poisson measure
$\mu$ on $[0, T] \times E$. We denote by $\mathbb{F}^W=(\mathcal{F}_t^W)_{t\leq T}$
(resp. $\mathbb{F}^\mu=(\mathcal{F}_t^\mu)_{t\leq T} $) the $\mathbb{P}$-augmentation of the natural filtration of $W$ (resp. $\mu$). As usual, we denote by $\mathcal{B}(X)$ the Borel set of topological set $X$. We introduce the following subset: $E^\eps :=\{e\in \R, \mbox{ s.t } / \ |e|\leq \eps\}$, $E_\eps :=\{e\in \R, \mbox{ s.t } / \ |e| > \eps\}$,  $E:=\R= E^\eps  \cup E_\eps$.
\2 \\
The martingale measure $\bar{\mu}$ is the compensated measure corresponding to Poisson random measure $\mu$, such that
$ \bar{\mu}(de,dr)=\mu(de,dr) - \nu(de)dr$, where $\nu$ is a Lévy measure on $E$ endowed with its Borel tribe $\mathcal{E}$. The Lévy measure $\nu$ will be assumed to satisfy $\nu(\R) =\infty$ and $\int_\R |e|^2\nu(de) < \infty $. Throughout this paper we deal with the measure $\bar{M}$ defined by
 \beqs \bar{M}(t, B )  = \int _{[0, t] \times B } e \bar{\mu}(dr,de), \s\s B \in \mathcal{B}(E) \enqs
which can  be considered as  a compensated Poisson random measure on $[0, T]\times E$ and   $\int _{[0, t] \times E} e \mu(dr,de)$  is a compound
Poisson random variable. We associate to $\bar{M}$ the $\sigma$-finite measure
\begin{equation}
m(B) := \int_B e^2 \nu(de) \s\s B \in \mathcal{B}(E).
\end{equation}
In particular, we have $\sigma(\eps)^2=m(E^\eps) $. \2 \\
The measure $\bar{M}$ is taken to drive the jump noise instead of $\bar{\mu}$, in the aim to adopt the concept of Malliavin calculus on the canonical
Lévy space from \cite{DI11}.
\2 \\
For some constant $K >0$, we consider four  $K$-Lipschitz functions with   bounded derivatives
$\beta: \mathbb{R} \rightarrow \mathbb{R}$, $b: \mathbb{R} \rightarrow \mathbb{R}$ , $g: \mathbb{R} \rightarrow \mathbb{R}$ and
$f :\Omega \times \mathbb{R} \times \mathbb{R} \times L^2(E, \mathcal{E}, \nu, \mathbb{R})\rightarrow \mathbb{R}$, where the first derivative of
$b$, $\beta$ and $g$ are bounded. \2 \\
Define $\rho$ to be a measurable function $ \rho: E \rightarrow \mathbb{R}$ such that:
        \begin{equation} \label{rho}
        \sup_{e \in E} \left|\rho(e) \right|< K.
        \end{equation}
For any $p\geq 2$ we consider the following class of  processes:
\begin{itemize}
  \item $S^p$ is the set of real valued adapted rcll process $Y$ such that:
        $$ \|Y\|_{S^p} := \mathbb{E}\left(\sup_{0\leq t\leq T}|Y_t|^p\right)^\frac{1}{p} < \infty. $$
  \item $H^p$ is the set of progressively measurable $\mathbb{R}$-valued processes $Z$ such that:
        $$ \|Z\|_{H^p} := \left(\mathbb{E}\left(\int_0^T|Z_r|^2dr\right)^{\frac{p}{2}}\right)^\frac{1}{p} < \infty. $$
  \item $L^p$ is the set of $\mathcal{P}\otimes \mathcal{E}$ measurable map $U: \Omega \times [0, T] \times E \rightarrow \mathbb{R} $
     such that:
        $$ \|U\|_{L^p} := \left(\mathbb{E}\int_0^T \int_E |U_r(e)|^p \nu(de)dr\right)^\frac{1}{p} < \infty. $$
  \item The space $\mathcal{B}^p := \mathcal{S}^p \times H^p \times L^p$ is endowed with the norm
        $$ \|(Y, Z, U)  \|_{\mathcal{B}^p} := \Big( \| Y \|_{S^p}^p + \| Z\|_{H^p} ^p + \| U \|_{L^p}^p \Big)^\frac{^1}{p}.$$
  \item $M^{2,p}$ the class of square integrable random variable $F$ of the form:
    $$ F= \mathbb{E}F + \int_0^TU_r dW_r +\int_0^T \int_E \psi(r,e)\bar{\mu}(de, dr),$$
    where $u$ (resp. $\psi$) is a progressively measurable (resp. measurable) process satisfying
    $\sup_{t \leq T} \mathbb{E}|u_t|^p < \infty$ (resp. $\sup_{t \leq T} \mathbb{E}\int_E|\psi(t,e)|^p\nu(de) < \infty$).
 \end{itemize}
\2
\subsection{Approximation scheme }

\vspace{2mm}
In this subsection, we show that the approximation error
\beqs Err^2_\eps(Y, V)&:=& \mathbb{E}\left[\sup_{t\leq T}|Y_t  - Y_t^\eps|^2 \right] + \mathbb{E} \left[\sup_{t\leq T} \left| \int_0^tV_r dR_r - \int_0^tZ^\eps_rdW_r \right|^2 \right] \\
    && + \mathbb{E}\left[\int_0^T\int_{E_\eps} |V_r - U_r^\eps (e)|^2m(de)dr \right],
\enqs
converges to $0$ as $\eps$ goes to $0$.
\begin {Theorem} \label{ThXepsToX}
Under the space  $(\Omega, \mathcal{F}, \mathbb{P}) $,
\begin{enumerate}
  \item There exist a solution $X$ on $[0, T]$ of
  \beq X_t &=& X_0+ \int_0^t b(X_r)dr + \int_0^t \beta(X_{r^-})dL_r, \enq
    where $X_0 \in \R$.
  \item There exist a solution $X^\eps$ on $[0, T]$ of
\beq X_t^\eps&=& X_0^\eps+ \int_0^t b(X_r^\eps)dr + \int_0^t \beta(X_r^\eps)\sigma(\eps)dW_r + \int_0^t \int_{E_\eps}\beta(X_{r^-}^\eps) \bar{M}(dr,de), \enq
    where $X^\eps_0 \in \R$.
\end{enumerate}
Moreover \beq             \mathbb{E} \left(\sup_{0 \leq t \leq T} | X_t - X_t^\eps |^2 \right) &\leq& C \sigma(\eps)^2 \label{X-Xeps}. \enq

\end {Theorem}
For the proof, we state firstly the following Lemma
\begin{Lemma}
\label{AprioriEstimate}
On the space $(\Omega, \mathcal{F}, \mathbb{P}) $, fixing  $\eps > 0$, we have for $p \geq 2$:
    \beq
                \mathbb{E} \left(\sup_{0 \leq t \leq T}  | X_t|^p \right)          &<& \infty.\\
                \mathbb{E} \left(\sup_{0 \leq t \leq T}  | X_t^\eps|^p \right) &<& \infty .
    \enq
\end{Lemma}
{\bf Proof.} We denote by $C$ a constant whose value may change from line to line. Using Jensen's inequality, Burkholder-Davis-Gundy inequality and Lipschitz property of $b$ and $\beta$ we have:
\beqs
 \mathbb{E} \sup_{0 \leq s \leq t} | X_s|^p &\leq& C \sup_{0 \leq s \leq t} \left( \mathbb{E} \left[ |X_0|^p + \left(\int_0^s  b(X_r)dr \right)^p + \int_0^s \int_E \beta(X_r)\bar{M}(dr,de) \right]^p \right) \\
 &\leq& C \left( |X_0|^p + \int_0^t \mathbb{E}\Big[ |b(X_0)| + |X_r| + |X_0|  \Big]^p dr\right. \\
        && + \left. \int_0^t \int_E \mathbb{E}\Big[ |\beta(X_0)| + |X_0| + |X_r| \Big]^p e^p \nu(de)dr\right) \\
 &\leq& C \left( |X_0|^p + |b(X_0)|^p+ |\beta(X_0)|^p+ \int_0^t \mathbb{E}\Big[ \sup_{0 \leq u \leq r}|X_u|^p dr\Big] \right),
 \enqs
where $C$ depends on $t$, $b(X_0)$ and $\beta(X_0)$.
We conclude the first assertion by Gronwall's Lemma. \2 \\ Following the same arguments, we obtain the second assertion. \ep
\2 \\
{\bf Proof of Theorem \ref{ThXepsToX}}  The existence and uniqueness of such SDEs was already studied in the literature see e.g. \cite{FK89} and \cite{A04}. Then it remains to prove the estimate (\ref{X-Xeps}). \2 \\
Using  Jensen's inequality leads to :
\beqs
 \mathbb{E} \sup_{0 \leq u \leq t} | X_u - X^\eps_u|^2 &\leq& C \left(\mathbb{E} \left[  \int_0^t |b(X_r) - b(X_r^\eps)|dr \right]^2 \right.\\
                &&+ \mathbb{E} \left[  \int_0^t \int_{E_\eps} |\beta(X_r) - \beta(X_r^\eps) | \bar{M}(dr,de) \right]^2\\
                && \left.+\mathbb{E} \left[  \left(\int_0^t \int_{E^\eps} |\beta(X_r)|\bar{M}(de, dr) \right)^2 + \left( \int_0^t |\beta(X_r^\eps) \sigma(\eps)|dW_r \right)^2  \right]\right).
 \enqs
By Burkholder-Davis-Gundy inequality, we get
 \beqs
  \mathbb{E} \sup_{0 \leq u \leq t} | X_u - X^\eps_u|^2  &\leq& C \left(t\mathbb{E} \left[ \int_0^t \left(b(X_r) - b(X_r^\eps) \right)^2dr \right] \right.\\
                &&+ \mathbb{E} \left[ \int_0^t \int_{E_\eps} \left(\beta(X_r) - \beta(X_r^\eps) \right)^2m(de)dr\right] \\
                &&+  \left.\mathbb{E} \left[  \int_0^t \int_{E^\eps} \beta(X_r)^2m(de) dr \right] + C\mathbb{E} \left[ \int_0^t \beta(X_r^\eps)^2 \sigma(\eps)^2dr \right] \right).
\enqs
By  Lipschitz property of $b$ and $\beta$
\beqs
 \mathbb{E} \sup_{0 \leq u \leq t} | X_u - X^\eps_u|^2&\leq& C  \left(\left[\int_0^t \mathbb{E} \left(X_r - X_r^\eps \right)^2 dr \right]  \right.\\
            && \left.+\sigma(\eps)^2 \mathbb{E}\int_0^t \left[ \beta^2(X_0)+\beta^2(X^\eps_0)+ |X^\eps_r|^2 +|X_r|^2\right]dr \right)\\
            &\leq& C \left[\int_0^t \mathbb{E} \left[\sup_{u \leq r} |\delta X_u|^2 \right] dr + \sigma(\eps)^2 \right].
\enqs\\
The result follows from Gronwall's Lemma. \ep
\2 \\
Finally,  we can now state  the main result of this section.
\begin{Theorem} \label{ThErrEps}
Under the space  $(\Omega, \mathcal{F}, \mathbb{P}) $,

\begin{enumerate}
  \item There exist a unique pair  $(Y, V) \in \mathcal{S}^2 \times H^2$, which solves the BSDE:
\begin{eqnarray} \label{Y^eps}
Y_t &=& g(X_T) + \int_t^Tf(\Theta_r)dr -\int_t^T  V_rdL_r ,
\end{eqnarray}
where $\Theta :=\Big(X, Y, \int_{E}\rho(e)Ve \nu(de)\Big)$.
  \item For a fixed $\eps >0$, There exist a unique solution $(Y^\eps, Z^\eps, U^\eps) \in \mathcal{B}^2$ of  the following BSDE:
\begin{eqnarray}\label{Y}
Y_t^\eps &=& g(X^\eps_T) + \int_t^Tf (\Theta^\eps_r)dr -\int_t^T  Z^\eps_rdW_s  - \int_t^T \int_{E_\eps} U_r^\eps  (e) \bar{M}(dr,de).
\end{eqnarray}
with $\Theta^\eps :=\Big(X^\eps, Y^\eps, \Gamma^\eps \Big)$ and $ \Gamma^\eps := \int_{E_\eps}\rho(e)U^\eps(e)e \nu(de)$.\\

\end{enumerate}
Moreover, if \hspace{1mm} $\sup_{t \leq T} \mathbb{E}|V_t|^2 < \infty$, then there exist a constant $C$ such that:
\beq\label{errorYEpsl}
Err^2_\eps(Y, V)\leq C \sigma(\eps)^2.
\enq
\end{Theorem}
\begin{Remark}\label{RemYeps}
Observe that
\beqs \mathbb{E}\left[\sup_{t\leq T}|Y_t  - Y_t^\eps|^2 \right] &+& \mathbb{E} \left[\sup_{t\leq T} \left| \int_t^T \int_\mathbb{R}V_r \bar{M}(de,dr)  \right.\right.\\
&&  \s\s - \left.\left.\int_0^tZ^\eps_rdW_r  - \int_0^T\int_{E_\eps}U^\eps_r(e) \bar{M}(de,dr) \right|^2\right]\\
 &\leq& Err^2_\eps(Y, V)\\
    &\leq& \sigma(\eps)^2.
\enqs
Which shows clearly the convergence of the approximated scheme (\ref{FBSDE_Eps}) to (\ref{FBSDE_L}).
\end{Remark}
{\bf Proof. of Theorem \ref{ThErrEps}} Existence and uniqueness of the solutions of BSDEs (\ref{Y^eps}) and (\ref{Y})  was already proved, see e.g  \cite{BBP97}.\2 \\
We are going to prove inequality (\ref{errorYEpsl}). By Itô's formula applied to $|\delta Y|^2 := |Y - Y^\eps|^2$ yields :
\begin{eqnarray}
\label{Ito1}
|\delta Y_t |^2 &+& \int_t^T {Z_r^\eps}^2dr + \int_t^T \int_{E_\eps} \left(U_r^\eps (e) - V_r\right)^2m(de)dr \nonumber\\
    &=& |g(X_T) - g(X_T^\eps)|^2 + \sigma(\eps)^2 \int_t^T V_r^2dr - \int_t^T \delta Y_rZ^\eps_r dW_r \nonumber\\
    &&+ 2\int_t^T\delta Y_r \left(f\left(\Theta_r \right) - f\left(\Theta_r^\eps\right)\right)dr + 2\int_t^T\int_{E^\eps}\Big[(\delta Y_{s^-} + V_r)^2 - \delta Y_{s^-}^2\Big]\bar{M}(de, dr) \nonumber\\
    &&- 2\int_t^T\int_{E_\eps}\Big[(\delta Y_{s^-} + U_r^\eps (e) - V_r)^2 - \delta Y_{s^-}^2\Big]\bar{M}(de, dr).
\end{eqnarray}
Taking expectation in both hand-side of the above equality we get
\begin{eqnarray}
\label{Ito2}
&\mathbb{E}& \left[\delta Y_t ^2 + \int_t^T {Z_r^\eps}^2dr + \int_t^T \int_{E_\eps} \left(U_r^\eps (e) - V_r\right)^2m(de)dr \right] \nonumber\\
    &&= \mathbb{E}\left[|g(X_T) - g(X_T^\eps)|^2 + \sigma(\eps)^2 \int_t^T V_r^2dr + 2\int_t^T\delta Y_r \left(f\left(\Theta_r \right) - f\left(\Theta_r^\eps\right)\right)dr \right].\nonumber
\end{eqnarray}
From Lemma \ref{AprioriEstimate}, Lipschitz property of $g$ and Jensen inequality  we obtain:
\beqs
&\mathbb{E}& \left[|Y_t - Y_t^\eps|^2 + \int_t^T {Z^\eps}_r^2dr + \int_t^T \int_{E_\eps} \left(U_r^\eps (e) - V_r\right)^2m(de)dr \right] \\
    &\leq& C \left(\sigma(\eps)^2 + K\mathbb{E} \int_t^T(Y_r - Y_r^\eps)(|X_r - X_r^\eps| +|Y_r - Y_r^\eps|) dr \right.\\
    && + K\mathbb{E} \int_t^T\Big[(Y_r - Y_r^\eps)\int_{E_\eps} \rho(e)e|U_r^\eps (e) - V_r|\nu(de)\Big]dr \\
    &&+ \left.K\mathbb{E} \int_t^T\Big[(Y_r - Y_r^\eps)\int_{E^\eps } \rho(e)V_re\nu(de)\Big]dr \right).
\enqs
Using the fact that $ab \leq \alpha a^2 + \frac{1}{\alpha}b^2$ for some $\alpha > 0$, yields to

    \beqs&\mathbb{E}& \left[|Y_t - Y_t^\eps|^2 + \int_t^T {Z^\eps}_r^2dr + \int_t^T \int_{E_\eps} \left(U_r^\eps (e) - V_r\right)^2m(de)dr \right] \\
&\leq& C \left(\sigma(\eps)^2 + K(1 + \alpha^2 +\gamma^2 + \eta^2)\mathbb{E}\int_t^T(Y_r - Y_r^\eps)^2 dr  + \frac{K}{\alpha^2}\mathbb{E}\int_t^T|X_r - X_r^\eps|^2dr \right.\\
    && + \left.\frac{K}{\gamma^2}\mathbb{E}\int_t^T\int_{E_\eps} \rho(e)^2|U_r^\eps (e) - V_r|^2m(de)dr + \frac{K}{\eta^2}\mathbb{E}\int_t^T\int_{0 \leq |e| \leq \eps } \rho(e)^2V_r^2m(de)dr\right),
\enqs
\\
where $\alpha$ and $\gamma$ are two constants taken such that $\frac{K}{\alpha^2} = \frac{K^3}{\gamma^2}=\frac{1}{2}$, we then get

\begin{eqnarray}
\mathbb{E} \Big[|Y_t - Y_t^\eps|^2 &+& \int_t^T {Z^\eps}_r^2dr + \int_t^T \int_{E_\eps} \left(U_r^\eps (e) - V_r\right)^2m(de)dr \Big] \nonumber\\
    &\leq& C\left( \sigma(\eps)^2 + (K + 2K^2 + 2K K^2) \mathbb{E}\int_t^T(Y_r - Y_r^\eps)^2 dr \right).
\end{eqnarray}
\\
Using Gronwall's Lemma, we deduce that
\begin{eqnarray}
\label{Estimate1}
\mathbb{E} |Y_t - Y_t^\eps|^2 \leq C \sigma(\eps)^2.
\end{eqnarray}
Plugging this estimate in the previous upper bound, we get
\begin{eqnarray}
\label{Estimate2}
\mathbb{E}\int_0^T {Z_r^\eps}^2dr + \mathbb{E}\int_0^T \int_{E_\eps} \left(U_r^\eps (e) - V_r\right)^2m(de)dr \leq C \sigma(\eps)^2,
\end{eqnarray}
Then
\begin{eqnarray}
\mathbb{E} |\delta Y_t |^2 + \mathbb{E}\int_0^T {Z_r^\eps}^2dr + \mathbb{E}\int_0^T \int_{E_\eps} \left(U_r^\eps (e) - V_r\right)^2m(de)dr \leq C \sigma(\eps)^2.
\end{eqnarray}
Now using Burkholder-Davis-Gundy inequality, we have
\begin{eqnarray}
\label{Estimate3}
\mathbb{E} \sup_{t \leq T}|\delta Y_t |^2 + \mathbb{E}\int_0^T {Z_r^\eps}^2dr + \mathbb{E}\int_0^T \int_{E_\eps} \left(U_r^\eps (e) - V_r\right)^2m(de)dr \leq C \sigma(\eps)^2.
\end{eqnarray}
From other side, it follows by Burkholder-Davis-Gundy inequality and (\ref{Estimate3}) that:
\begin{eqnarray}
\label{Estimate4}
\mathbb{E}\left[\sup_{t \leq T} \left| \int_0^t Z^\eps_rdW_r - \int_0^t V_rdR_r  \right|^2 \right]  &\leq&  C \mathbb{E} \left[\sup_{t \leq T} \left| \int_0^t Z^\eps_rdW_r \right|^2  + \sup_{t \leq T} \left| \int_0^t V_rdR_r  \right|^2 \right] \nonumber \\
            &\leq& C \mathbb{E} \left[\int_0^T {Z_r^\eps}^2dr  +  \sigma(\eps^2)\int_0^T V_r^2dr\right] \nonumber \\
            &\leq& \mathbb{E} \int_0^T {Z_r^\eps}^2dr  +  C\sigma(\eps^2).
\end{eqnarray}
The result now follows by combining (\ref{Estimate3}) and (\ref{Estimate4}).
\ep\\
\begin{Remark}
One can show that $\sup_{t \leq T} \mathbb{E}|V_t|^2 < \infty$, following the same arguments as in Remark 2.8 in  \cite{HNS11}.\end{Remark}

\begin{Remark}
In the BSDE (\ref{Y}), each time we change $\eps$, there exist a unique pair $(Z^\eps, U^\eps)$ of predictable process, such that the BSDE (\ref{Y})
has a solution.
\end{Remark}

\section{Forward-backward Euler scheme}
\setcounter{equation}{0} \setcounter{Assumption}{0}
\setcounter{Theorem}{0} \setcounter{Proposition}{0}
\setcounter{Corollary}{0} \setcounter{Lemma}{0}
\setcounter{Definition}{0} \setcounter{Remark}{0}

In this section, we discretize the solution $(X^\eps, Y^\eps, Z^\eps, \Gamma^\eps)$ of (\ref{FBSDE_Eps}) by $({X}^\pi, {Y}^\pi, {Z}^\pi, {\Gamma}^\pi)$ defined by induction in (\ref{Approximation}) and then we show the convergence of $({X}^\pi, {Y}^\pi, {Z}^\pi, {\Gamma}^\pi)$ to the solution of (\ref{FBSDE_R}).
Thus let us recall some definition and notation.\\
For each $t\in [t_i, t_{i+1})$, we define:
\begin{eqnarray}
\bar{Z}_t=n \mathbb{E} \left[ \int_{t_i}^{t_{i+1}} Z_s ds / \mathcal{F}_{t_i}   \right], \s \s
\bar{\Gamma}_t=n \mathbb{E} \left[ \int_{t_i}^{t_{i+1}} \Gamma_sds / \mathcal{F}_{t_i}   \right],
\end{eqnarray}
and
\begin{eqnarray}
\bar{Z}^\pi_{t_i} =n \mathbb{E} \left[ \int_{t_i}^{t_{i+1}} Z_s^\pi ds / \mathcal{F}_{t_i}   \right], \s \s
\bar{\Gamma}^\pi_{t_i} =n \mathbb{E} \left[ \int_{t_i}^{t_{i+1}} \Gamma_s^\pi ds / \mathcal{F}_{t_i}   \right].
\end{eqnarray}
The process $\bar{Z}_{t_i}$ and $\bar{\Gamma}_{t_i}$  (resp. $\bar{Z}^\pi_{t_i}$ and $\bar{\Gamma}^\pi_{t_i}$)
can be interpreted as the best approximation of $Z_{t_i}$ and $\Gamma_{t_i}$  (resp. $Z^\pi_{t_i}$ and $\Gamma^\pi_{t_i}$).
We know from Bouchard and Elie \cite{BE08}, that FBSDE (\ref{Y^eps}) has a backward Euler scheme taking the form:

\begin{equation}
\label{Approximation}
\left\{
\begin{array}{lll}
\bar{X}^\pi_{t_i+1} &=& \bar{X}^\pi_{t_i}+\frac{1}{n}b(\bar{X}^\pi_{t_i}) + \sigma(\eps)\Delta W_{i+1} + \int_{E_\eps}\beta(\bar{X}^\pi_{t_i}) \bar{M}(de,(t_i, t_{i+1}] ) \\
\bar{Z}^\pi_{t}&=&n\mathbb{E} \left[ \bar{Y}^\pi_{t_{i+1}} \Delta W_{i+1} / \mathcal{F}_{t_i} \right] \\
\bar{\Gamma}^\pi_{t}&=&n\mathbb{E} \left[ \bar{Y}^\pi_{t_{i+1}} \int_{E_\eps} \rho(e) \bar{M} (de,(t_i, t_{i+1}] )   / \mathcal{F}_{t_i} \right] \\
\bar{Y}^\pi_{t}&=&\mathbb{E} \left[ \bar{Y}^\pi_{t_{i+1}}  / \mathcal{F}_{t_i} \right] + \frac{1}{n} f(\bar{X}^\pi_{t_i}, \bar{Y}^\pi_{t_i}, \bar{\Gamma}^\pi_{t_i} )
\end{array}
\right.
\end{equation} \\
for which the discretization error:
\begin{eqnarray} \label{errorYPhi}
\overline{Err}_n(Y^\eps, Z^\eps, \Gamma^\eps) &:=& \left\{ \sup_{t\leq T} \mathbb{E} \left[ |Y^\eps_t - \bar{Y}^\pi_t|^2 \right] + \|Z^\eps - \bar{Z}^\pi\|^2 _{H^2} +  \|\Gamma^\eps - \bar{\Gamma}^\pi \|^2 _{H^2}  \right\}^\frac{1}{2} \nonumber \\
                & \leq& C n^{-1/2},
\end{eqnarray}
 converges to $0$ as the discretization step $\frac{T}{n}$ tends to $0$. Means that the discretization scheme (\ref{Approximation}) achieves the optimal convergence rate $n^{-1/2}$. The regularity of $Z^\eps$ and $\Gamma^\eps$ has been studied in $L^2$ sense in \cite{BE08} when the terminal value is a functional of forward diffusion.\vspace{2mm} \\
It is well known also that
\beqs
\max_{i<n}\E \left[\sup_{t\in[t_i, t_{i+1}]}|X^\eps_{t_i} - \bar{X}_{t_i}^\pi|^2 \right] \leq Cn^{-1/2}.
\enqs
Our aim in this part, is to show that the approximation-discretization error between BSDE (\ref{FBSDE_R}) and (\ref{Approximation}):
\begin{eqnarray}
\label{DefErrorEpsN}
\overline{Err}^2_{(n,\eps)}(Y, V) &:=& \sup_{t\leq T} \mathbb{E} \left[ |Y_t- \bar{Y}^\pi_t|^2 \right] + \sup_{t \leq T} \mathbb{E} \left|\int_0^t V_rdR_r - \int_0^t\bar{Z}_r^\pi dW_r \right|^2 \nonumber \\
                    && +  \|\Gamma - \bar{\Gamma}^\pi \|^2 _{H^2},
\end{eqnarray}
converges to $0$ as $(\eps, n) \rightarrow (0, \infty)$. \2 \\
The first main result of this paper is:
\begin{Proposition} \label{PropErrNEps}
Assuming Lipschitz property of coefficients $b$ and $\beta$, the approximation-discretization error defined in (\ref{DefErrorEpsN}) is bounded by:
    \begin{equation}\label{ErrNEps}
    \overline{Err}_{(n,\eps)}(Y, V) \leq C \left(n^{-1/2} + \sigma(\eps) \right).
    \end{equation}
Means that:
\begin{equation}
\overline{Err}_{(n,\eps)}(Y, V) \underset{(n, \eps)\rightarrow (\infty, 0)}{ \longrightarrow }0. \\
\nn
\end{equation}

\end{Proposition}
{\bf Proof.} From (\ref{DefErrorEpsN}), Jensen inequality and Burkholder-Davis-Gundy inequality, we have:
\begin{eqnarray}
\overline{Err}^2_{(n,\eps)}(Y, V) &\leq& C \Big(\sup_{t\leq T} \mathbb{E} \left[ |Y_t - Y^\eps_t |^2 +  |Y^\eps_t -\bar{Y}^\pi_t|^2 \right]  +  \|\Gamma     - \Gamma^\eps \|^2 _{H^2} + \|\Gamma^\eps - \bar{\Gamma}^\pi \|^2 _{H^2}  \nn \\
        && \hspace{5mm} + \sup_{t \leq T} \mathbb{E} \left|\int_0^t V_rdR_r - \int_0^t Z^\eps_r dW_r \right|^2 + \|Z^\eps -  \bar{Z}^\pi \|^2_{H^2}\Big).
\end{eqnarray}
From other side, it follows from Hölder inequality that:
\begin{eqnarray}
\int_{E_\eps}\rho(e)e(V_r- U_r^\eps (e))\nu(de) \leq K \left(\int_{E_\eps}(V_r- U_r^\eps (e))^2m(de)\right)^\frac{^1}{2} \nu\Big({E_\eps}\Big)^\frac{^1}{2}.
\end{eqnarray}
Recalling that $\nu(E_\eps) < \infty$  so that $\nu$ has a.s. only a finite number of big jumps on $[0,T]$. Combining the two last inequalities with (\ref{errorYPhi}) leads to:
\beqs
\overline{Err}^2_{(n,\eps)}(Y, V) &\leq& C\left(n^{-1} + \sup_{t\leq T} \mathbb{E} \left[ |Y_t - Y^\eps_t |^2 \right] \right.+ \sup_{t \leq T} \mathbb{E} \left|\int_0^t V_rdR_r - \int_0^t Z_r dW_r \right|^2 \\
        && \hspace{10mm} +  \int_0^T\int_{E^\eps} \E|V_r |^2 m(de) dr+ \left.\mathbb{E}\left[\int_0^T\int_{E_\eps} |U_r^\eps (e) - V_r|^2m(de)dr \right] \right).
\enqs
By Theorem \ref{ThErrEps} we get:
\begin{eqnarray}
\overline{Err}_{(n,\eps)}(Y, V) &\leq& C \left(n^{-1/2} + \sigma(\eps) \right),
\end{eqnarray}
where $C$ depends on $K$. \ep

\vspace{2mm}
\begin{Remark}
In the general case, as we neglect the small jump, the Brownian part in (\ref{FBSDE_Eps}) disappears. In this case the assertion (\ref{ErrNEps}) can be replaced by:
\beqs
    \overline{Err}_{(n,\eps)}(Y, V) \leq C n^{-1/2}.
\enqs
\end{Remark}
\begin{Remark} Taking $\eps = n^{-1/2}$, we obtain the optimal convergence rate $n^{-1/2}$  in (\ref{ErrNEps}),
\beqs    \overline{Err}_{(n,\eps)}(Y, V) \leq C n^{-1/2},
\enqs
which is exactly the approximation error  in \cite{BE08}.
\end{Remark}
\section{A discrete scheme via Malliavin derivatives}
\setcounter{equation}{0} \setcounter{Assumption}{0}
\setcounter{Theorem}{0} \setcounter{Proposition}{0}
\setcounter{Corollary}{0} \setcounter{Lemma}{0}
\setcounter{Definition}{0} \setcounter{Remark}{0}
In this section, we generalize the new discrete scheme recently introduced by Hu, Nualard and Song \cite{HNS11} from a general BSDE, to our framework of decoupled
forward-backward SDEs with jumps. For this aim, we use the Malliavin derivatives of $Y$ to derive the discrete scheme.
We first fix a regular grid $ \pi:= \{t_i:=iT/n, i=0, ..., n \}$ on $[0,T] $ and approximate the forward SDE $X^\eps$ in (\ref{FBSDE_Eps}) by its Euler scheme $X^\pi$ already defined in (\ref{Approximation}).
It'is hard to prove existence and convergence of Malliavin derivatives of $\bar{X}^\pi$.
However, to avoid this problem, we can instead consider the continuous-time version of the Euler scheme, then we define the function $\phi$ for each $t\in [0, T]$:
\begin{equation}
\phi^n_t := \max \{t_i, \  i=0, ..., n. \ /  \ t_i \leq t\},
\end{equation}
for which we associate:
  \beqs
  X^\pi_t&:=& X^\pi_{\phi^n_t} + \sigma(\eps)b(X^\pi_{\phi^n_t})( t -\phi^n_t)+  \sigma(\eps)\beta(X^\pi_{\phi^n_t})(W_t - W_{\phi^n_t}) + \int_{E_\eps} \beta(X^\pi_{\phi^n_t})\bar{M}(de,(t, \phi^n_t]).
  \enqs
It could be written as
  \beqs
  X^\pi_t&:=& X_0 +\int_0^t b(X^\pi_{\phi^n_r})dr + \int_0^t \sigma(\eps)\beta(X^\pi_{\phi^n_r})dW_r+ \int_0^t\int_{E_\eps} \beta(X^\pi_{\phi^n_r})\bar{M}(de,dr).
  \enqs
It is well known that under Lipschitz property of the coefficients
\begin{eqnarray}\label{XXpi}
\mathbb{E}\left[\sup_{t\in [0, T]}|X^\eps_t - X^\pi_t|^p\right]^{1/p} \leq C \pi^{1/2}.
\end{eqnarray}
The Malliavin derivatives of the continuous-time version of  Euler scheme for $\theta \leq s$ a.e. are:
\beqs
D_\theta X^\pi_t &=& \int_\theta^t \partial_xb(X^\pi_{\phi^n_r})D_\theta X^\pi_{\phi^n_r}dr + \int_\theta^t \int_{E_\eps} \partial_x \beta(X^\pi_{\phi^n_r})D_\theta X^\pi_{\phi^n_r}\bar{M}(de,dr)+ \sigma(\eps)\beta(X^\pi_\theta)\\
                    && + \sigma(\eps)\int_\theta^t \partial_x\beta(X^\pi_{\phi^n_r})D_\theta X^\pi_{\phi^n_r}dW_r,\\
D_{\theta,e}X^\pi_t &=& \int_\theta^t D_{\theta,e}b(X^\pi_r)dr + \int_\theta^t \int_{E_\eps} D_{\theta,e}\beta(X^\pi_r)\bar{M}(de,dr)+ \sigma(\eps)\int_\theta^t D_{\theta,e}\beta(X^\pi_r)dW_r \\
                    && +\beta(X^\pi_\theta).
\enqs
We introduce some additional assumptions:
    \begin{description}
          \item[(A1)] $f(t,y,\gamma)$ doesn't depend on $x$.
          \item[(A2)] The first derivative of $b$ and $\beta$ and $g$ is a K-Lipschitz function
            \beqs
            |b'(x) - b'(y)| + |\beta'(x) - \beta'(y)| +|g'(x) - g'(y)|  \leq K|x-y|.
            \enqs
          \item[(A3)] $f(t,y,\gamma)$ is linear with respect to $t$, $y$ and $\gamma$. Moreover, there exist three bounded functions $f_1$, $f_2$ and $f_3$ such that :
          \begin{eqnarray}
          f(t, y, u)= f_1(t) + f_2(t)y + f_3(t)\gamma.
          \end{eqnarray}

        \end{description}
\begin{Lemma} \label{LemDXpi}
Under Lipschitz continuity of $b$ and $\beta$, we have for any $q \geq 1$
\begin{eqnarray}
\sup_{0\leq \theta \leq T}\sup_{n \geq 1 } \mathbb{E}\left[ \sup_{\theta \leq t \leq T}\| D_\theta X^\pi_t\|^{2q} \right] &<& \infty, \\
\sup_{0\leq \theta \leq T}\sup_{n \geq 1 } \mathbb{E}\left[ \sup_{\theta \leq t \leq T}\| D_{\theta,e} X^\pi_t\|^{2q} \right]&<&  \infty.
\end{eqnarray}
\end{Lemma}
For the proof see the Appendix. \2 \\
We then derive the following theorem
\begin{Theorem} Under assumption (A2), Lipschitz continuity of $b$ and $\beta$ and  for any $p \geq 2$, we have,
\begin{eqnarray}
\mathbb{E}\left[\sup_{t\in [0, T]}|D_\theta X^\eps_t - D_\theta X^\pi_t|^p\right]^{1/p} &\leq& C_p \pi^{1/2} \label{DXXpi}.\\
\mathbb{E}\left[\sup_{t\in [0, T]}|D_{\theta,e}X_t^\eps - D_{\theta,e} X^\pi_t|^p\right]^{1/p} &\leq& C_p \pi^{1/2}.
\end{eqnarray}
\end{Theorem}

{\bf Proof.} Using Burkholder-Davis-Gundy inequality, Jensen inequality, inequality (\ref{XXpi}) and  Lemma \ref{LemDXpi}
\beqs
\mathbb{E}\left[\sup_{s\in [0, t]}|D_\theta X_t^\eps - D_\theta X^\pi_t|^p\right]  & \leq& C_p\mathbb{E}\left[\int_0^t \Big|\partial_xb(X^\pi_{\phi^n_r})D_\theta X^\pi_{\phi^n_r} - \partial_xb(X_r)D_\theta X^\eps_r \Big|^p dr \right.\\
            && \s \s + \int_0^t \int_{E_\eps} \Big|\partial_x \beta(X^\pi_{\phi^n_r})D_\theta X^\pi_{\phi^n_r} - \partial_x \beta(X^\eps_r)D_\theta X^\eps_r\Big|^p\nu(de)dr\\
            && \s \s + \sigma^p(\eps)\left(\int_0^t \Big|\partial_x\beta(X^\pi_{\phi^n_r})D_\theta X^\pi_{\phi^n_r}- \partial_x\beta(X^\eps_r)D_\theta X^\eps_r\Big|^2dr\right)^{p/2} \\
            && \left.\s \s + \sigma^p(\eps)\Big|\beta(X^\pi_\theta) -\beta(X^\eps_\theta)\Big|^p \right].
\enqs
which leads to
\beqs
& \leq& C_p\mathbb{E}\left[\int_0^t \big|D_\theta X^\eps_r\big|^p \Big[|\partial_xb(X^\eps_r)-\partial_xb(X^\pi_{\phi^n_r})|^p +|\partial_x\beta(X^\eps_r)-\partial_x\beta(X^\pi_{\phi^n_r})|^p \Big]dr \right.\\
            && \s \s   + \int_0^t\big|D_\theta X^\eps_r - D_\theta X^\pi_{\phi^n_r}\big|^p \Big[|\partial_xb(X^\pi_{\phi^n_r})|^p  + |\partial_x\beta(X^\pi_{\phi^n_r})|^p\Big] dr \\
            && \s \s  + |X^\pi_\theta - X^\eps_\theta|^p \Big]\\
& \leq& C_p\mathbb{E}\left(\pi^{p/2} +  \int_0^t \sup_{u\in [0, r]}|D_\theta X^\eps_u - D_\theta X^\pi_u|^p dr \right).
\enqs
We conclude by using Gronwall's Lemma. Following the same arguments, we prove the second assertion.\\ \ep \2 \\
Now we derive the discrete scheme using the expression of $Z^\eps$ and $U^\eps$ as the trace of Malliavin derivatives of $Y$. From equation (\ref{D_the0Y}) and (\ref{D_theeY}), the two Malliavin derivatives $D_\theta Y_t$ , $D_{\theta,e} Y_t$ could be expressed as:
\begin{eqnarray}
D_\theta Y^\eps_t    &=& \mathbb{E}\left(\mathcal{E}_{t,T}\partial_x g(X_T^\eps)D_\theta X_T^\eps + \int_t^T \mathcal{E}_{t,r}\partial_xf(\Theta_r^\eps)D_\theta  X_r^\eps dr / \mathcal{F}_t \right)  \label{Zt}\\
D_{\theta,e} Y^\eps_t&=& \mathbb{E}\left(\mathcal{E}^e_{t,T}D_{\theta,e} g(X_T^\eps) + \int_t^T \mathcal{E}^e_{t,r}\alpha_{\theta,r}D_{\theta,e} X_r^\eps dr / \mathcal{F}_t \right)\label{Ute},
\end{eqnarray}
where
\beqs
\mathcal{E}_{t,r}   &=& exp\left\{\int_t^r\left(\partial_yf(\Theta^\eps_u) - \frac{1}{2}\int_{E_\eps}\partial_\gamma f^2(\Theta^\eps_u)\rho^2(e)m(de) \right)du \right.\\
                    && \s \s + \left.\int_t^r\int_{E_\eps}\partial_\gamma f(\Theta^\eps_u)\rho(e)\bar{M}(de, du)  \right\} \\
\mathcal{E}^e_{t,r} &:=& exp\left\{\int_t^r \Big[\alpha_{\theta,u} - \frac{1}{2}\alpha^2_{\theta,u}\int_{E_\eps}\rho^2(e)m(de) \Big]du + \int_t^r\int_{E_\eps}\alpha_{\theta,u}\rho(e)\bar{M}(de, du)  \right\}
\enqs
and
\beqs
\alpha_{\theta,r}  &:=&\frac{f(\Theta^\eps_r + D_{\theta,e}\Theta^\eps_r) - f(\Theta^\eps_r)}{D_{\theta,e}X_r^\eps + D_{\theta,e}Y_r+D_{\theta,e}\Gamma_r} 1_{\{D_{\theta,e}X_r^\eps + D_{\theta,e}Y_r + D_{\theta,e}\Gamma_r \neq 0\}}.
\enqs
\vspace{2mm}
Thus, we define our discrete scheme for $ i= n-1, ..., 1, 0.$ and $t\in [t_i, t_{i+1})$ by induction
\begin{equation}
\label{NewYPi}
\left\{
\begin{array}{lll}
Y^\pi_{t_i}&=&\mathbb{E}\Big[ Y^\pi_{t_{i+1}} + f(\Theta^\pi_{t_{i+1}})\Delta t_{i+1}/ \mathcal{F}_{t_i}\Big]\\
Z^\pi_{t_i}&=&\mathbb{E}\Big[\mathcal{E}^\pi_{t_{i+1}, t_n}\partial_x g(X^\pi_T) D_{t_i}X^\pi_T +\sum_{k=i}^{n-1}\mathcal{E}^\pi_{t_{i+1}, t_{k+1}} \partial_xf(\Theta^\pi_{t_{k+1}})D_{t_i}X^\pi_{t_{k+1}}\Delta t_{k}/ \mathcal{F}_{t_i}\Big]\\
\Gamma^\pi_{t_i}&=&\E\Big[\int_{E_\eps}\rho(e)\big[\mathcal{E}^{e,\pi}_{t_{i+1}, t_n}D_{t_i,e}g(X^\pi_T) +\sum_{k=i}^{n-1}\mathcal{E}^{e,\pi}_{t_{i+1}, t_{k+1}} \alpha^\pi_{t_i, t_{k+1}}D_{t_i,e}X^\pi_{t_{k+1}}\Delta t_k\big]\nu(de) / \mathcal{F}_{t_i}\Big]\\
\end{array}
\right.
\end{equation}
with terminal conditions \\
\s $Y^\pi_{t_n}=g(X^\pi_T)$,  \s $ Z^\pi_{t_n}=\sigma(\eps) \partial_x g(X^\pi_T)\beta(X^\pi_T) $,  \s $ U^\pi_{t_n,e}=g(X^\pi_T + \beta(X^\pi_T)) - g(X^\pi_T)$, \2 \\
where for any $0\leq i < j \leq n$,
\begin{eqnarray}
\mathcal{E}^\pi_{t_i, t_j} &=& exp\left\{\sum_{k=i}^{j-1}\int_{t_k}^{t_{k+1}}\Big[\partial_yf(\Theta^\pi_{t_k})  - \frac{1}{2}\int_{E_\eps}\partial_\gamma f^2(\Theta^\pi_{t_k})\rho^2(e)m(de) \Big]dr \right.\nonumber\\
                            && \hspace{10mm}  + \left.\sum_{k=i}^{j-1}\int_{t_k}^{t_{k+1}}\int_{E_\eps}\partial_\gamma f(\Theta^\pi_{t_k})\rho(e)\bar{M}(de, dr)  \right\} \label{Epii},\\
\mathcal{E}^{e,\pi}_{t_i, t_j}&=& exp\left\{\sum_{k=i}^{j-1}\int_{t_k}^{t_{k+1}}\Big[\alpha^\pi_{\theta,r,t_k} - \frac{1}{2}{\alpha^\pi}^2_{\theta,r,t_k}\int_{E_\eps}\rho^2(e)m(de)  \Big]dr\right.\nonumber\\
                            && \hspace{10mm} + \left.\sum_{k=i}^{j-1}\int_{t_k}^{t_{k+1}}\int_{E_\eps}\alpha^\pi_{\theta,r,t_k}\rho(e)\bar{M}(de, dr)  \right\}\label{Epiei},
\end{eqnarray}
and
\beqs
\alpha^\pi_{\theta,r,t_k}   &:=&\frac{f(\Theta^\pi_{t_k}+ D_{\theta,e}\Theta^\pi_{t_k}) - f(\Theta^\pi_{t_k})}{D_{\theta,e}X^\pi_{t_k} + D_{\theta,e}Y^\pi_{t_k}+D_{\theta,e}\Gamma^\pi_{t_k}} \times 1_{\{D_{\theta,e}X^\pi_{t_k} + D_{\theta,e}Y^\pi_{t_k} + D_{\theta,e}\Gamma^\pi_{t_k} \neq 0\}}, \\
\Gamma^\pi_{t_k}      &:=&\int_{E_\eps} U^\pi_{t_k,e}\rho(e)\nu(de),
\enqs
with $\Theta^\pi_{t_k} =\left(r,X^\pi_{t_k},Y^\pi_{t_k},\Gamma^\pi_{t_k} \right) $. \2 \\
We are going to compute the discretization error of our discrete scheme and prove the convergence of the above scheme. We recall the expression of the error between the solution of (\ref{FBSDE_Eps}) and (\ref{NewYPi}):
\beqs
Err^p_n(Y^\eps,Z^\eps,\Gamma^\eps)&:=& \mathbb{E} \max_{0 \leq i \leq n} \left[|Y^\eps_{t_i} -Y^\pi_{t_i}|^p + |Z^\eps_{t_i} -Z^\pi_{t_i}|^p + |\Gamma^\eps_{t_i} -\Gamma^\pi_{t_i}|^p  \right],
\enqs
where, $\Gamma^\pi_{t_i}=\int_{E_\eps} \rho(e)U^\pi_{t_i,e}\nu(de)$. \2 \\
We also recall the expression of discretization-approximation error between (\ref{FBSDE_R}) and (\ref{NewYPi})
\beqs
Err^2_{n, \eps}(Y, V) &:=& \max_{0 \leq i \leq n} \sup_{t\in[t_i, t_{i+1}]} \mathbb{E}  \left[|Y_t -Y^\pi_{t_i}|^2 \right] + \E\left|\int_0^TV_rdR_r  -\sum_{i=0}^{n-1}Z^\pi_{t_i}\Delta W_{t_i}\right|^2  \\
                                && + \sum_{i=0}^{n-1}\int_{t_i}^{t_{i+1}}\mathbb{E}|\Gamma_t -\Gamma^\pi_{t_i}|^2 dt.
\enqs
We conclude this section with the following Theorems whose proof are at the end of section 5.
\begin{Theorem} \label{ThErrorYi}
Under assumption \ref{asmpCoef}, we assume the existence of a  constant $L_3 >0 $ such that :
\begin{eqnarray}
|f(t_2,y,u) - f(t_1,y,u)  |\leq L_3|t_2 - t_1|^\frac{1}{2}.
\end{eqnarray}
Then there exist a positive constant $C$ independent of $\pi$ such that:
\beq \label{ErrorYi}
Err^p_n(Y^\eps,Z^\eps,\Gamma^\eps)\leq C_p |\pi|^{\frac{p}{2}-\frac{p}{2log\frac{1}{|\pi|}}}.
\enq
\end{Theorem}
\vspace{2mm}
The second main result of this paper is summarized in the following theorem:
\begin{Theorem}\label{ThErrorYEpsi}
Under the same assumptions as Theorem \ref{ThErrorYi}, we have
\begin{eqnarray}
Err_{n, \eps}(Y, V)\leq C \left( \sigma(\eps) +  |\pi|^{\frac{1}{2}-\frac{1}{2log\frac{1}{|\pi|}}}    \right).
\end{eqnarray}
\vspace{2mm}
\end{Theorem}
\begin{Remark} The importance of the above scheme is it can be adapted to a backward SDE  when the generator does not depends on the terminal value of a forward equation. \\
Consider the following backward stochastic differential equation driven by pure jump Lévy processes
 \beq \label{BSDExi} \hat{Y}_t&=& \hat{\xi} +\int_t^Tf\Big(r, \hat{Y}_r, \int_E \rho(e)\hat{V}_re \nu(de)\Big)dr - \int_t^T \hat{V}_r dL_r
 \enq
 Which we approximate by
\beq\label{BSDExieps} \hat{Y}_t^\eps &=& \hat{\xi} + \int_t^Tf\Big(r, \hat{Y}^\eps_r, \int_{E_\eps} \rho(e)\hat{V}_re \nu(de)\Big)dr -\int_t^T  \hat{Z}_r^\eps dW_s -\int_t^T \int_{E_\eps} \hat{U}_r^\eps  (e) \bar{M}(dr,de) \nn \\ \enq
We finally propose the below discrete time scheme, defined by terminal values $ \hat{Y}^\pi_{t_n}=\xi$, $\hat{Z}^\pi_{t_n}=D_T\xi$ and $ \hat{U}^\pi_{t_i}= D_{T,e}\xi$
\begin{equation}\label{BSDExipi}
\left\{
\begin{array}{lll}
\hat{Y}^\pi_{t_i}&=&\mathbb{E}\Big[ \hat{Y}^\pi_{t_{i+1}} + f(\hat{\Theta}^\pi_{t_{i+1}})\Delta t_{i+1}/ \mathcal{F}_{t_i}\Big]\\
\hat{Z}^\pi_{t_i}&=&\mathbb{E}\Big[ \mathcal{E}^\pi_{t_{i+1}, t_n}D_{t_i}\hat{\xi} +\sum_{k=i}^{n-1}\mathcal{E}^\pi_{t_{i+1}, t_{k+1}} \partial_xf(\hat{\Theta}^\pi_{t_{k+1}})D_{t_i}X^\pi_{t_{k+1}}\Delta t_{k}/ \mathcal{F}_{t_i}\Big]\\
\hat{\Gamma}^\pi_{t_i}&=&\E\Big[\int_{E_\eps}\rho(e)\big[\mathcal{E}^{e,\pi}_{t_{i+1}, t_n}D_{t_i,e}\hat{\xi } +\sum_{k=i}^{n-1}\mathcal{E}^{e,\pi}_{t_{i+1}, t_{k+1}} \alpha^\pi_{t_i, t_{k+1}}D_{t_i,e}X^\pi_{t_{k+1}}\Delta t_k\big]\nu(de)/ \mathcal{F}_{t_i}\Big]\\
\end{array}
\right.
\end{equation} \\
with $\hat{\Theta}^\pi_{t_k} =\left(r,\hat{Y}^\pi_{t_k},\hat{\Gamma}^\pi_{t_k} \right) $. \2 \\
Under the same assumptions of Theorem \ref{ThErrorYi} we prove the convergence of the system (\ref{BSDExipi}) to BSDE (\ref{BSDExi}). Moreover, we obtain the upper bound
 \beq
&& \max_{0 \leq i \leq n} \sup_{t\in[t_i, t_{i+1}]} \mathbb{E}  \left[|Y_t -Y^\pi_{t_i}|^2 \right] + \E\left|\int_0^TV_rdR_r  -\sum_{i=0}^{n-1}Z^\pi_{t_i}\Delta W_{t_i}\right|^2   + \sum_{i=0}^{n-1}\int_{t_i}^{t_{i+1}}\mathbb{E}|\Gamma_t -\Gamma^\pi_{t_i}|^2 dt \nn\\
                                &&\leq C \left( \sigma^2(\eps) +  |\pi|^{\frac{1}{}-\frac{1}{log\frac{1}{|\pi|}}}    \right).
                                \enq
\end{Remark}
\vspace{5mm}
\section{Malliavin calculus for FBSDEs}
\setcounter{equation}{0} \setcounter{Assumption}{0}
\setcounter{Theorem}{0} \setcounter{Proposition}{0}
\setcounter{Corollary}{0} \setcounter{Lemma}{0}
\setcounter{Definition}{0} \setcounter{Remark}{0}

For ease of notations, we shall denote throughout  this section  the process $(X^\eps, Y^\eps, Z^\eps, \Gamma^\eps)$ by $(X, Y, Z, \Gamma)$.\\
In this section, we study some regularity properties of the solution $(X, Y, Z, \Gamma)$. We recall the system (\ref{FBSDE_Eps}) using the new notations
\begin{equation} \label{FB}
\left\{
\begin{array}{lll}
X_t&=& X_0+ \int_0^t b(X_r)dr + \int_0^t \beta(X_r)\sigma(\eps)dW_r + \int_0^t \int_{E_\eps}\beta(X_{r^-}) \bar{M}(dr,de)\\
Y_t &=& g(X_T) + \int_t^Tf(\Theta_r)dr -\int_t^T  Z_r dW_s -\int_t^T \int_{E_\eps} U_r  (e) \bar{M}(dr,de)
\end{array}
\right.
\end{equation}
In fact, there are many methods to develop Malliavin calculus for Lévy processes. In our paper, we opt for the approach of Solé et al. \cite{SUV07}, based on a chaos decomposition in terms of multiple stochastic integrals with respect to the random measure $\bar{M}$. Adopting notation of \cite{DI11}, we will recall the suitable canonical space we adopt to our setting.
\vspace{2mm}\\
We start by introducing some additional notations and definitions. We assume that the probability space $(\Omega, \mathcal{F}, \mathbb{P})$ is the product of two canonical
spaces $(\Omega_W\times\Omega_\mu, \mathcal{F}_W\times\mathcal{F}_\mu,\mathbb{P}_W\times\mathbb{P}_\mu) $ and the filtration
$\mathbb{F}=(\mathcal{F}_t)_{t\in [0, T]}$ the canonical filtration completed for $\mathbb{P}$ (for details concerning this construction, see Section 2 in \cite{DI11}). \2 \\
We consider the finite measure $q$ defined on $[0, T]\times \mathbb{R}$ by
\beqs
q(B)=\int_{B(0)}dt + \int_{B'}e^2\nu(de)dt, \s\s B \in \mathcal{B}([0, T] \times \mathbb{R}).
\enqs
where $B(0)=\{ t\in [0, T]; (t, 0) \in B \}$, $B' = B - B(0) $ and the random measure $Q\in [0, T]\times \mathbb{R}$:
\beqs
Q(B)=\int_{B(0)}dW_t + \int_{B'}e\bar{\mu}(dt,de), \s\s B \in \mathcal{B}([0, T] \times \mathbb{R}).
\enqs
For $n\in \mathbb{N}$, a simple function $h_n =1_{E_1 \times ... \times E_n} $  with pairwise disjoints sets $E_1, ..., E_n \in \mathcal{B}([0, T] \times \mathbb{R}) $, we define:
\beqs
I_n(h_n)=\int_{([0, T]\times \mathbb{R})^n} h((t_1, e_1),...,(t_n, e_n)) Q(dt_1, de_1)\cdot...\cdot Q(dt_n, de_n).
\enqs \\
We Define the following spaces
\begin{enumerate}
  \item $\mathbb{L}^2_{T,q,n}(\mathbb{R})$ the space of product measurable deterministic functions $h: ([0, T]\times \mathbb{R})^n \rightarrow \mathbb{R}$ satisfying $\|h \|^2_{\mathbb{L}^2_{T,q,n}} < \infty$, where
    \beqs\|h \|^2_{\mathbb{L}^2_{T,q,n}}=:\int_{([0, T]\times \mathbb{R})^n} |h((t_1, e_1),...,(t_n, e_n))|^2 q(dt_1, de_1)\cdot...\cdot q(dt_n, de_n).\enqs
    \item $\mathbb{D}^{1,2}(\mathbb{R}) $ denote the space of $\mathbb{F}$-measurable random variables $H \in \mathbb{L}^2(\mathbb{R})$ with the representation $H=\sum_{n=0}^\infty I_n(h_n)$ and satisfying \beqs \sum_{n=0}^\infty nn!\|h_n \|^2_{\mathbb{L}^2_{T,q,n}} < \infty. \enqs
    \item $\mathbb{L}^{1,2}(\mathbb{R})$ denote the space of product measurable and $\mathbb{F}$-adapted processes $G:\Omega\times \mathbb{R}\rightarrow \mathbb{R}$ satisfying
      $$ \mathbb{E} \left(\int_{[0,T] \times \mathbb{R}} |G(s,y)|^2q(ds,dy)\right) < \infty$$
      $$ G(s,y) \in \mathbb{D}^{1,2}(\mathbb{R}), \mbox{ for } q-\mbox{a.e } (s, y) \in [0, T] \times \mathbb{R}$$
      $$ \mathbb{E} \left(\int_{([0,T] \times \mathbb{R})^2} |D_{t,z}G(s,y)|^2q(ds,dy)q(dt,dz)\right)< \infty.$$
      This space is endowed with the norm
      \beqs
      \|G \|^2_{\mathbb{L}^{1,q}} &=&\mathbb{E} \left(\int_{[0,T] \times \mathbb{R}} |G(s,y)|^2q(ds,dy) \right)\\
                                    &+& \mathbb{E} \left(\int_{([0,T] \times \mathbb{R})^2} |D_{t,z}G(s,y)|^2q(ds,dy)q(dt,dz)\right).
      \enqs
\end{enumerate}
We should mention that the derivative $D_{t,0}$ coincide  with $D_t$ the classical Malliavin derivative with respect to Brownian motion. \vspace{2mm}\\
To study the regularity of $Z$ and $U$, we shall also introduce the following assumption:
\begin{Assumption}\label{asmpCoef}
For $2 \leq p \leq \frac{q}{2} $
\begin{enumerate}
    \item The generator $f$ has continuous and uniformly bounded first and second order partial derivative with respect to $x$, $y$ and $\gamma$.
    \item For each $(x, y, \gamma) \in \mathbb{R}^3$, $\partial_xf(\Theta), \partial_yf(\Theta)$ and $\partial_\gamma f(\Theta)$
  belong to $\mathbb{L}^{1,2}$ and satisfy
  \begin{eqnarray}
  \sup_{0 \leq \theta \leq T} \mathbb{E} \left( \int_\theta^T |D_\theta \partial_if(\Theta_r)|^2dr\right)^\frac{q}{2} &<& \infty, \\
  \sup_{0 \leq \theta \leq T}\sup_{0 \leq u \leq T} \mathbb{E} \left( \int_{\theta\vee u}^T |D_u D_\theta \partial_if(\Theta_r)dr|^2 \right)^{\frac{q}{2}}&<& \infty,
  \end{eqnarray}
    where $i:=x, y, \gamma.$\\
    There exist a constant $K>0$ such that for any $e\in \big(\mathbb{R}- \{0\}\big)$, $t\in [0, T]$ and \\$0 \leq \theta, u \leq t \leq T $:
    \begin{eqnarray}
    \mathbb{E}  |D_\theta g(X_T)- D_u g(X_T)|^p &\leq& K |\theta - u|^\frac{p}{2} \label{AsmpDtheg}.\\
    \mathbb{E} |D_{\theta,e} g(X_T)- D_{u,e} g(X_T)|^p &\leq& K |\theta - u|^\frac{p}{2} \label{asmpDtheeg}.\\
    \mathbb{E} \left( \int_t^T |D_\theta f(\Theta_r)- D_u f(\Theta_r)|^2\right)^\frac{p}{2} &\leq& K |\theta - u|^\frac{p}{2}. \label{AsmpDthef}\\
    \mathbb{E} \left( \int_t^T |D_{\theta,e} f(\Theta_r)- D_{u,e} f(\Theta_r)|^2\right)^\frac{p}{2} &\leq& K |\theta - u|^\frac{p}{2}. \label{asmpDtheef}
    \end{eqnarray}
\end{enumerate}
\end{Assumption}
Additionally to Assumption \ref{asmpCoef}, we assume that
\begin{Assumption}
\label{AsmpAlpha}
 For any $\lambda >0$ and $q \geq 1$, we consider three progressive measurable processes $\{\alpha_t\}_{0 \leq t \leq T}$, $\{\beta_t\}_{0 \leq t \leq T}$
  and $\{\gamma_t\}_{0 \leq t \leq T}$ such that:
            \beqs \mathbb{E} \ exp \left( \lambda \int_0^T\Big(|\beta_r| + \gamma_r^2\big)dr \right) &<& \infty,     \\
             \sup_{0 \leq t \leq T} \mathbb{E} \Big(|\alpha_t|^q + | \gamma_t|^q\Big) &<& \infty. \enqs
 \end{Assumption}
\begin{Proposition}\label{propE}
Under Assumption \ref{AsmpAlpha}, the discontinuous semi-martingale  ${\mathcal{E}_t}$:
\begin{equation}
d\mathcal{E}_t=\mathcal{E}_t \beta_tdt + \mathcal{E}_t\gamma_t\int_{E_\eps}\rho(e)\bar{M}(de, dt),
\end{equation}
has the following properties
\begin{enumerate}
  \item $\mathbb{E}\underset{0 \leq t \leq T} {\sup}\mathcal{E}_t^n < \infty$, for any  $ n\in \mathbb{R}$.
  \item The process $ \mathcal{Z}_t :=  \mathcal{E}_t^{-1} $ satisfies the following linear SDE:
    \begin{equation*}
    \frac{d\mathcal{Z}_t}{\mathcal{Z}_t} = \left(-\beta_t + \gamma_t^2\int_{E_\eps} \rho^2(e)m(de)\right)dt - \gamma_t\int_{E_\eps}\rho(e)\bar{M}(de,dt).
    \end{equation*}
  Moreover, we have for any $p \geq 2$:
    \begin{equation}
    \mathbb{E} | \mathcal{Z}_t - \mathcal{Z}_s |^p \leq C |t-s|^p.
    \end{equation}
\end{enumerate}
\end{Proposition}
{\bf Proof.} $\mathcal{E}$ could be written as :
\begin{equation*}
\mathcal{E}_t = exp\left\{\int_0^t\Big[\beta_r - \frac{1}{2}\int_{E_\eps}\gamma_r^2\rho^2(e)m(de) \Big]dr + \int_0^t\int_{E_\eps}\gamma_r\rho(e)\bar{M}(de, dr)  \right\}.
\end{equation*}
Under Assumption \ref{AsmpAlpha}, we get the first assertion. The second assertion is deduced from first one, Hölder inequality and
Burkholder-Davis-Gundy inequality.
\vspace{5mm}
\ep\\
The following Theorem constitutes the main tool to prove Theorem \ref{ThTrace}.
\begin{Theorem}\label{linear BSDEs}
Suppose that $\mathcal{E}_T X_T$  and $ \int_0^T\alpha_r D_\theta X_rdr$ are in $M^{2,q} $. The following linear BSDE
    \begin{eqnarray}
       Y_t&=&g(X_T)X_T+ \int_t^T  \left[ \alpha_r X_r + \beta_r Y_r + \gamma_r\Gamma_r \right]dr - \int_t^TZ_rdW_r \nonumber  \\
            && -\int_t^T\int_{E_\eps} U_r(e)\bar{M}(dr,de), \s  0 \leq t \leq T
    \end{eqnarray}
    has a unique solution $(Y, Z, U)$ and there is a constant $C>0$ such that
    \begin{equation}\label{Yt-Ys}
    \mathbb{E} |Y_t-Y_s|^p \leq C|t-s|^\frac{p}{2} \s \mbox{for all } s,t\in [0, T].
    \end{equation}
\end{Theorem}
{\bf Proof.} Applying Itô's formula to $\mathcal{E}_tY_t$, we obtain
\beqs
d(\mathcal{E}_tY_t) &=&  -\mathcal{E}_t\alpha_t X_T  dt + \mathcal{E}_tZ_tdW_t  + \mathcal{E}_t \int_{E_\eps} (Y_t\gamma_t \rho(e) + U_t(e) )\bar{M}(de,dt). \enqs
Then
\beqs
Y_t= \mathbb{E}\left(\mathcal{E}_{t,T}g(X_T)X_T + \int_t^T \mathcal{E}_{t,r}\alpha_r X_r  dr / \mathcal{F}_t \right),
\enqs
where $\mathcal{E}_{t,r} = \mathcal{Z}_t\mathcal{E}_r$. \\
For $0 \leq s \leq t \leq T $, we have:
    \begin{eqnarray}
    \mathbb{E} |Y_t-Y_s|^p &\leq& 3^{p-1} \mathbb{E} \left| \mathbb{E}\left( \mathcal{E}_{t,T}g(X_T) X_T/\mathcal{F}_t \right) - \mathbb{E} \left( \mathcal{E}_{s,T}g(X_T)X_T /\mathcal{F}_s \right)\right|^p \nonumber  \\
                           &&+ 3^{p-1} \mathbb{E} \left| \mathbb{E}\left( \int_t^T\mathcal{E}_{t,r}\alpha_r X_r dr /\mathcal{F}_t \right) - \mathbb{E} \left( \int_s^T\mathcal{E}_{s,r}\alpha_r X_r dr /\mathcal{F}_s\right)\right|^p \nonumber  \\
                           &=&  3^{p-1} (I_1 + I_2).
    \end{eqnarray}
By adapting the argument of Theorem 2.3 in \cite{HNS11} and recall Remark \ref{RemXTheta}, we can immediately show that $I_1  \leq C |t-s|^{\frac{p}{2}} $ and $I_2 \leq C |t-s|^{\frac{p}{2}}$. \vspace{5mm}\ep


\subsection{Malliavin calculus on the Forward SDE}
In this section, we recall well-known properties on forwards SDEs, concerning the Malliavin derivatives of the solution of a forward SDE with jump, stated in Nualart \cite{N95} in the case of SDE without jumps and in Petrou \cite{P08} in case of a Lévy process. The following theorem can be found in \cite{P08}.
\begin{Theorem} $ $\\
Let $X$ be the solution of forward SDE (\ref{FB}). Then, for all $t\in [0, T]$ and $(\theta, e) \in [0, T]\times (\mathbb{R} \backslash \{0\})$,  the Malliavin derivatives of $X$ satisfy
\begin{eqnarray}
    D_{\theta}X_t &=& \int_\theta^t \partial_xb(X_r)D_{\theta}X_rdr + \int_\theta^t \int_{E_\eps}\partial_x \beta(X_{r})D_{\theta}X_r\bar{M}(dr,de)\nonumber\\
                    &&+  \sigma(\eps)\beta(X_\theta)+ \sigma(\eps)\int_\theta^t \partial_x\beta(X_r)D_{\theta}X_rdW_r, \s 0 \leq \theta \leq t \leq T.
\end{eqnarray}
And
\begin{eqnarray}
    D_{\theta,e}X_t &=& \int_\theta^t D_{\theta,e}b(X_r)dr + \int_\theta^t \int_{E_\eps}D_{\theta,e}\beta(X_{r})\bar{M}(dr,de)\nonumber\\
                    &&+  \beta(X_\theta)+ \sigma(\eps)\int_\theta^t D_{\theta,e}\beta(X_r)dW_r, \s 0 \leq \theta \leq t \leq T.
\end{eqnarray}
For all $\theta >t$, we have $D_{\theta,e}X_t  = D_{\theta}X_t=0 $ a.s.
\end{Theorem}
\begin{Remark} \label{RemXTheta}
Using standard arguments as in the proof of Lemma \ref{LemDXpi}, we can prove the following a priori estimate:
\beqs
    \sup_{0\leq \theta \leq T} \mathbb{E} \left[\sup_{0\leq t \leq T} D_\theta X_t\right] &<& \infty,\\
    \sup_{0\leq \theta \leq T , e\in \mathbb{R}^*}\mathbb{E} \left[\sup_{0\leq t \leq T} D_{\theta,e}X_t \right] &<& \infty,\\
    \sup_{0\leq u \leq T}\sup_{0\leq \theta \leq T} \mathbb{E} \left[\sup_{0\leq t \leq T} D_u D_\theta X_t\right] &<& \infty,\\
    \sup_{0\leq u \leq T}\sup_{0\leq \theta \leq T, e\in \mathbb{R}^*} \mathbb{E} \left[\sup_{0\leq t \leq T} D_u D_{\theta,e} X_t\right] &<& \infty.\\
\enqs
\end{Remark}
\subsection{Malliavin calculus on the Backward SDE}
In this section, we recall some result of Malliavin derivatives applied to BSDE especially established in \cite{EPQ97} and \cite{DI11} in the aim to generalize the result of Theorem 2.6 in \cite{HNS11}.
\begin{Theorem} \label{ThTrace} Assume that assumption \ref{asmpCoef} hold. There exist a unique solution \\$\{(Y_t,Z_t, U_t(e))\}_{0 \leq t \leq T, e\in (\mathbb{R} - \{0\})} $ of BSDE (\ref{FB}), such that:
\begin{enumerate}
  \item The first version of Malliavin derivative $\{(D_{\theta}Y_t, D_{\theta}Z_t, D_{\theta}U_t(e))_{0 \leq \theta, t \leq T, e\in (\mathbb{R} - \{0\})}$ of the solution $\{(Y_t,Z_t, U_t(e))\}_{0 \leq t \leq T, e\in (\mathbb{R} - \{0\})} $ satisfies the following linear BSDE :
    \begin{eqnarray}\label{D_the0Y}
    D_{\theta}Y_t &=& \partial_xg(X_T)D_{\theta}X_T + \Int_t^T f^{\theta}(X_r, Y_r, \Gamma_r) dr- \int_t^T  D_{\theta}Z_rdW_s     \\
     && -\int_t^T \int_{E_\eps} D_{\theta}U_r(e) \bar{M}(dr,de), \s 0 \leq \theta \leq t \leq T \nonumber
    \end{eqnarray}\\
    where $f^{\theta}(\Theta):= \partial_x f(\Theta) D_{\theta}X_r +\partial_y f(\Theta) D_{\theta}Y_r + \partial_\gamma f(\Theta) D_{\theta}\Gamma_r.$ \2 \\
     Moreover $(D_tY_t)_{0\leq t \leq T}$ is a version of  $(Z_t)_{0\leq t \leq T}$ :
     \begin{equation}\label{DYversion}
     Z_t = D_tY_t. \s \mbox{a.s.}
     \end{equation}
 \item The second version of Malliavin derivative $(D_{\theta,z}Y_t, D_{\theta,z}Z_t, D_{\theta,z}U_t(e))_{0 \leq \theta, t \leq T, (e,z)\in (\R - \{0\})^2}$
 of the solution $(Y_t,Z_t, U_t(z))_{0 \leq t \leq T, z\in (\R - \{0\})} $ satisfies the following linear BSDE
     \begin{eqnarray}\label{D_theeY}
    D_{\theta,z}Y_t &=& g(X_T + D_{\theta,z}X_T) - g(X_T) + \int_t^T [f(\Theta_r + D_{\theta,z}\Theta_r) - f(\Theta_r)]dr \\
                    && - \int_t^T  D_{\theta,z}Z_rdW_s  -\int_t^T \int_{E_\eps} D_{\theta,z}U_r (e) \bar{M}(dr,de), \s 0 \leq \theta \leq t \leq T. \nn
    \end{eqnarray}
    Moreover $(D_{t,e}Y_t)_{0\leq t \leq T, e\in (\R - \{0\})}$ is a version of  $(U_t(z))_{0\leq t \leq T, z\in (\R - \{0\})}$:
    \begin{equation}\label{DeYversion}
     U_t(z) = D_{t,e}Y_t \s \mbox{a.s.}
     \end{equation}

    And for $(\theta, e) \in [0, T]\times \mathbb{R}$
    \beqs
     D_{\theta,z}Y_t = D_{\theta,z}Z_t = D_{\theta,z} U_t(z) = 0, \s  0 \leq t < \theta \s (e,z) \in \R \times (\R - \{0\}).
    \enqs

 \item There exist a constant $C >0$ such that for all $s, t \in [0, T]$:
          \begin{eqnarray}
          \mathbb{E}|Z_t - Z_s|^p &\leq& C|t-s|^\frac{p}{2} \label{Zt-Zs}\\
          \mathbb{E}|\Gamma_t - \Gamma_s|^p &\leq& C |t-s|^\frac{^p}{2} \label{Gt-Gs}.
          \end{eqnarray}
\end{enumerate}
\2
\end{Theorem}

{\bf Proof.} Existence and uniqueness of solution is similar to Proposition 5.3 in \cite{EPQ97} and Theorem 4.1 in \cite{DI11}. Then we focus our attention
to prove inequalities (\ref{Zt-Zs}) and (\ref{Gt-Gs}). \\

$\blacksquare$ We first prove that $\mathbb{E}|Z_t- Z_s|^p\leq C |t-s|^\frac{^p}{2}$. \2 \\
Let $C>0$ be a constant  independent of $s$ and $t$, whose value vary from line to line. From (\ref{DYversion}) we have:
     \beqs      Z_t- Z_s= D_tY_t - D_sY_s. \enqs
Then
 \beqs      \mathbb{E}|Z_t- Z_s|^p\leq \mathbb{E} |D_tY_t - D_sY_t|^p + \mathbb{E} |D_sY_t - D_sY_s|^p .\enqs \\
 \noindent {\bf Step 1}: Estimate $\mathbb{E}|D_tY_t - D_sY_t|^p$. \2 \\
     From Lemma (\ref{LemAprioriEstimateBSDEs}), equation (\ref{D_the0Y}) and assumption (\ref{AsmpDtheg})-(\ref{AsmpDthef}), we obtain
       \begin{eqnarray}
       \mathbb{E}|D_tY_t - D_sY_t|^p &+& \mathbb{E} \left(\int_t^T|D_tZ_r - D_sZ_r|^2dr\right)^\frac{p}{2} \nonumber\\
       &+&  \mathbb{E} \left( \int_t^T \int_{E_\eps} |D_tU_r (e) - D_sU_r (e)|^2 m(de)dr \right)^{\frac{p}{2}} \nonumber \\
       &\leq& C \mathbb{E} \Big( |D_tg(X_T) - D_sg(X_T)|^p  \Big) \nonumber\\
       &&    + C\mathbb{E}\left( \int_t^T|D_tf(r,X_r, Y_r, \Gamma_r)-D_sf(r, X_r, Y_r, \Gamma_r)|^2dr\right)^\frac{p}{2} \nonumber\\
       &\leq &C|t-s|^\frac{p}{2} \label{step1Z}.
       \end{eqnarray} \\
\noindent {\bf Step 2}: Estimate $\mathbb{E}|D_sY_t - D_sY_s|^p$. \vspace{2mm}
\\
We recall the expression of $\mathcal{E}_t$
\begin{equation}
\mathcal{E}_t = exp\left\{\int_0^t\Big[\beta_r - \frac{1}{2}\int_{E_\eps}\gamma_r^2\rho^2(e)m(de) \Big]dr + \int_0^t\int_{E_\eps}\gamma_r\rho(e)\bar{M}(de, dr)  \right\}.
\end{equation}
Denote $\beta_r = \partial_yf(\Theta_r)$ and $\gamma_r = \partial_\gamma f(\Theta_r) $. For any $ 0\leq \theta \leq t \leq T$, we have:
\beqs
D_\theta \mathcal{E}_t = &\mathcal{E}_t& \left\{\int_\theta^t\int_{E_\eps} \rho(e)\Big[\partial_{\gamma x}f(\Theta_r)D_\theta X_r + \partial_{\gamma y}f(\Theta_r)D_\theta Y_r   + \partial_{\gamma \gamma}f(\Theta_r)D_\theta \Gamma_r\Big] \bar{M}(de, dr)\right.\\
&+& \int_\theta^t\int_{E_\eps} \left[ \partial_{xy}f(\Theta_r) -\rho^2(e) \gamma_r \partial_{x\gamma} f(\Theta_r)\right] D_\theta X_r m(de)dr  \\
&+& \int_\theta^t\int_{E_\eps} \left[ \partial_{yy}f(\Theta_r) -\rho^2(e) \gamma_r \partial_{y\gamma} f(\Theta_r)\right] D_\theta Y_r m(de)dr\\
&+&\left.\int_\theta^t\int_{E_\eps} \left[ \partial_{\gamma y}f(\Theta_r) -\rho^2(e) \gamma_r \partial_{\gamma \gamma} f(\Theta_r)\right] D_\theta \Gamma_r m(de)dr \right\}.
\enqs
From other side, by induction on chain rule:
\begin{equation}\label{f2}
D_{\theta,e}f^2(\Theta)=f^2(\Theta+D_{\theta,e}\Theta)-f^2(\Theta).
\end{equation}
Using this in the previous equality
\beqs
D_{\theta,e} \mathcal{E}_t &=& \mathcal{E}_t \left(exp\left\{\int_\theta^t \Big[\partial_yf(\Theta_r+D_{\theta,e}\Theta_r )-\partial_yf(\Theta_r)\right.\right.\\
                && \s \s - \frac{1}{2}\int_{E_\eps}\rho^2(e) \Big[(\partial_\gamma f(\Theta_r + D_{\theta,e}\Theta_r))^2-\gamma_r^2\Big]m(de) \Big]dr  \\ &&\s\s+\gamma_\theta\rho(e)+ \left.\left.\int_\theta^t\int_{E_\eps}\rho(e)\Big[\partial_\gamma f(\Theta_r + D_{\theta,e}\Theta_r ) - \gamma_r\Big]\bar{M}(de, dr) \right\} -1\right).
\enqs
From Proposition \ref{propE}, Assumption \ref{asmpCoef}, Hölder inequality and Burkholder-Davis-Gundy inequality, we can show for any $p < q $ that:
\begin{eqnarray}
\sup_{\theta \in [0,T], e \in \mathbb{R}}\mathbb{E}\sup_{\theta \leq t\leq T}|D_{\theta,e}\mathcal{E}_t |^p &<& \infty. \label{DeE}
\end{eqnarray}
Now, by Clark-Ocone formula (See Es-Sebaiy and Tudor \cite{ET08}) on $\mathcal{E}_TD_\theta X_T$:
\beqs
\mathcal{E}_TD_\theta X_T   &=& \mathbb{E}(\mathcal{E}_TD_\theta X_T) + \int_0^T\mathbb{E}\Big(D_r(\mathcal{E}_TD_\theta X_T)/\mathcal{F}_r\Big)dW_r \\
                            && + \int_0^T\int_{E_\eps}\mathbb{E}\Big(D_{r,e}(\mathcal{E}_TD_\theta X_T)/\mathcal{F}_r\Big)\bar{M}(dr,de) \\
                            &=& \mathbb{E}(\mathcal{E}_TD_\theta X_T) + \int_0^T u^\theta_rdW_r + \int_0^T\int_{E_\eps} v^\theta_{r,e}\bar{M}(dr,de).
\enqs
Where
\beqs
u^\theta_r&:=& \mathbb{E}\Big(D_r\mathcal{E}_TD_\theta X_T + \mathcal{E}_TD_rD_\theta X_T/\mathcal{F}_r\Big) \\
v^\theta_{r,e}&:=& \mathbb{E}\Big( D_{r,e}\mathcal{E}_TD_\theta X_T +\mathcal{E}_TD_{r,e}D_\theta X_T + D_{r,e}\mathcal{E}_TD_{r,e}D_\theta X_T /\mathcal{F}_r\Big).
\enqs
Thus it remains to prove that
\beqs \sup_{\theta \in [0,T]}\sup_{r \in [0,T]}|u^\theta_r|^p &<& \infty \\
\sup_{\theta \in [0,T]}\sup_{r\in [0,T], e\in \mathbb{R}}|v^\theta_{r,e}|^p &<& \infty. \enqs
By Hölder inequality
\beqs
\mathbb{E}|v^\theta_{r,e}|^p &=& \mathbb{E} \Big|\mathbb{E}\Big(D_{r,e}\mathcal{E}_TD_\theta X_T +\mathcal{E}_TD_{r,e}D_\theta X_T + D_{r,e}\mathcal{E}_TD_{r,e}D_\theta X_T /\mathcal{F}_r\Big)\Big|^p \\
        &\leq& 3^{p-1}\Big(\mathbb{E} |D_{r,e}\mathcal{E}_TD_\theta X_T|^p + \mathbb{E}|\mathcal{E}_TD_{r,e}D_\theta X_T|^p + \mathbb{E}|D_{r,e}\mathcal{E}_TD_{r,e}D_\theta X_T |^p \Big) \\
        &\leq& 3^{p-1}\left( \left(\mathbb{E} |D_{r,e}\mathcal{E}_T|^{\frac{pq}{q-p}} \right)^\frac{q-p}{q} \Big(\mathbb{E} |D_\theta X_T|^q \Big)^\frac{p}{q}  + \left(\mathbb{E}|\mathcal{E}_T|^{\frac{pq}{q-p}}\right)^\frac{q-p}{q}\Big(\mathbb{E}|D_{r,e}D_\theta X_T|^q \Big)^\frac{p}{q} \right.  \\
        && \hspace{10mm} + \left.\left(\mathbb{E}|D_{r,e}\mathcal{E}_T|^{\frac{pq}{q-p}}\right)^\frac{q-p}{q}\Big(\mathbb{E}|D_{r,e}D_\theta X_T |^q\Big)^\frac{p}{q}\right).
\enqs
Combining  (\ref{DeE}) and Remark \ref{RemXTheta}, we deduce that $\sup_{\theta \in [0,T]}\sup_{r\in [0,T], e\in \mathbb{R}}|v^\theta_{r,e}|^p < \infty$. Following the same arguments we conclude that $\sup_{\theta \in [0,T]}\sup_{r \in [0,T]}|u^\theta_r|^p < \infty$. As consequence $\mathcal{E}_TD_\theta X_T $ belongs to $ M^{2,p}$. Therefore, by Theorem \ref{linear BSDEs} we conclude
     \begin{equation} \label{step2Z}
        \mathbb{E}|D_sY_t - D_sY_s|^p \leq C |t-s|^\frac{^p}{2}.
     \end{equation}
Now, combining (\ref{step1Z}) and (\ref{step2Z}), we finally obtain for some constant $C>0$
    \begin{equation}
      \mathbb{E}|Z_t- Z_s|^p= \mathbb{E}|D_sY_t - D_sY_s|^p \leq C |t-s|^\frac{^p}{2}.
     \end{equation} \\ \vspace{2mm}
$\blacksquare$ We prove that $\mathbb{E}|\Gamma_t - \Gamma_s|^p \leq C |t-s|^\frac{^p}{2}$.\\
By Hölder inequality
\beqs
\mathbb{E}|\Gamma_t - \Gamma_s|^p   &=& \mathbb{E}\left| \int_{E_\eps} \rho(e)(U_t(e) - U_s(e))\nu(de)\right|^p \\
                                    &\leq& \mathbb{E} \left(\left(\int_{E_\eps} \left| U_t(e) - U_s(e)\right|^p\nu(de)\right) \left(\int_{E_\eps} |\rho(e)|^\frac{p}{p-1}\nu(de) \right)^{p-1}\right) \\
                                    &\leq& C \int_{E_\eps} \mathbb{E} \left| U_t(e) - U_s(e)\right|^p\nu(de) \\
                                    &=&C \int_{E_\eps} \mathbb{E} \left| D_{t,e}Y_t - D_{s,e}Y_s\right|^p\nu(de) \\
                                    &^\leq&C \int_{E_\eps} \mathbb{E} \Big[\left| D_{t,e}Y_t - D_{s,e}Y_t\right|^p  + \left| D_{s,e}Y_t - D_{s,e}Y_s\right|^p \Big]\nu(de).
\enqs
\vspace{2mm}
\noindent {\bf Step 3.}: We prove that $\mathbb{E}\left| D_{t,e}Y_t - D_{s,e}Y_t\right|^p \leq C |t-s|^\frac{^p}{2}$. \\
Under assumption (\ref{asmpDtheeg}), (\ref{asmpDtheef}) and from Lemma \ref{LemAprioriEstimateBSDEs}

       \begin{eqnarray}
       \mathbb{E}|D_{t,e}Y_t - D_{s,e}Y_t|^p &+& \mathbb{E} \left(\int_t^T|D_{t,e}Z_r - D_{s,e}Z_r|^2\right)^\frac{p}{2} \nonumber\\
       &&  +  \mathbb{E} \left( \int_t^T \int_{E_\eps} |D_{t,e}U_r (e) - D_{s,e}U_r (e)|^2 \nu(de)dt \right)^{\frac{p}{2}} \nonumber \\
       &\leq & C \mathbb{E} \left[ |D_{t,e}\xi - D_{s,e}\xi|^p \right]\nonumber \\
       &&   + C\mathbb{E}\left( \int_t^T|D_{t,e}f(r, Y_r, U_r )-D_{s,e}f(r, Y_r, U_r)|^2dr\right)^\frac{p}{2} \nonumber\\
       &\leq & C|t-s|^\frac{p}{2}. \label{step1G}
       \end{eqnarray}\\
\vspace{2mm}
{\bf Step 4.}: We prove that $\mathbb{E}\left| D_{s,e}Y_t - D_{s,e}Y_s\right|^p\leq C |t-s|^\frac{^p}{2}$. \\
We can write BSDE (\ref{D_theeY}) as:
\begin{eqnarray} \label{linear BSDEsDThete}
    D_{\theta,e}Y_t &=& G(X_T)D_{\theta,e}X_T + \int_t^T \alpha_{\theta,r} \Big[D_{\theta,e}X_r +  D_{\theta,e} Y_r+ D_{\theta,e}\Gamma_r\Big]dr - \int_t^T  D_{\theta,e}Z_rdW_s  \nn \\
    &&-\int_t^T \int_{E_\eps} D_{\theta,e}U_r  (e) \bar{M}(dr,de),
          \end{eqnarray}
where
\beqs
G(X_T) &:=& \frac{g(X_T + D_{\theta,e}X_T) - g(X_T)}{D_{\theta,e}X_T}1_{\{D_{\theta,e}X_T \neq 0 \}} \\
\alpha_{\theta,r}  &:=&\frac{f(\Theta_r + D_{\theta,e}\Theta) - f(\Theta_r)}{D_{\theta,e}X_r + D_{\theta,e}Y_r+D_{\theta,e}\Gamma_r} 1_{\{D_{\theta,e}X_r + D_{\theta,e}Y_r + D_{\theta,e}\Gamma_r \neq 0\}}.
\enqs \\
Then, from Lipschitz continuity of $f$, we have $ \sup_{0 \leq t \leq T}\mathbb{E}|\alpha_{\theta,t}|^p < \infty $. It remain to show that $ \mathcal{E}_TD_{\theta,e}X_T$ belongs to $M^{2,p} $.
In fact, by Clark-Ocone formula applied to $ \mathcal{E}_TD_{\theta,e}X_T$:
\beqs
\mathcal{E}_TD_{\theta,e}X_T&=& \mathbb{E}(\mathcal{E}_TD_{\theta,e}) + \int_0^T\mathbb{E}\Big(D_r(\mathcal{E}_TD_{\theta,e} X_T)/\mathcal{F}_r\Big)dW_r \\
                            && + \int_0^T\int_{E_\eps}\mathbb{E}\Big(D_{r,e}(\mathcal{E}_TD_{\theta,e}X_T)/\mathcal{F}_r\Big)\bar{M}(dr,de) \\
                            &=& \mathbb{E}(\mathcal{E}_TD_{\theta,e}X_T) + \int_0^T \tilde{u}^\theta_rdW_r + \int_0^T\int_{E_\eps} \tilde{v}^\theta_{r,e}\bar{M}(dr,de),
\enqs
with
\beqs
\tilde{u}^\theta_r&:=& \mathbb{E}\Big(D_r\mathcal{E}_TD_{\theta,e}X_T + \mathcal{E}_TD_rD_{\theta,e}  X_T /\mathcal{F}_r\Big) \\
\tilde{v}^\theta_{r,e}&:=& \mathbb{E}\Big(D_{r,e}\mathcal{E}_TD_{\theta,e}X_T +\mathcal{E}_TD_{r,e}D_{\theta,e}X_T + D_{r,e}\mathcal{E}_TD_{r,e}D_{\theta,e}X_T /\mathcal{F}_r\Big).
\enqs
Following the same argument as \emph{Step 2} and using Remark (\ref{RemXTheta}) we prove that  $$\sup_{\theta \in [0,T], e\in \mathbb{R}} \mathbb{E} \left(|\tilde{u}^\theta_r|^p + |\tilde{v}^\theta_{r,e}|^p \right) < \infty.$$
Therefore, $ \mathcal{E}_TD_{\theta,e}X_T$ belongs to $M^{2,p} $. Finally, we apply once  again  the result of Theorem \ref{linear BSDEs} to BSDE (\ref{linear BSDEsDThete}) we get:
\begin{equation}
\mathbb{E}|D_{s,e}Y_t - D_{s,e}Y_s | \leq C |t-s|^\frac{p}{2}.
\end{equation}
The result then follows.\ep
\vspace{5mm}\\
We now complete the proof of Section 3 \\
{\bf Proof. of Theorem \ref{ThErrorYi}} We adapt the proof of Theorem 5.2 in \cite{HNS11}. \\
Let $i=n-1, ..., 1, 0 $.\2 \\
{\bf Step 1. }: We show that $\mathbb{E} \left[\sup_{0 \leq i \leq n}|\delta Z^\pi_{t_i}|^p \right] \leq C |\pi|^{p-1}$. \\
Denote \beqs \delta Z^\pi_{t_i} = Z_{t_i} -Z^\pi_{t_i}. \enqs
Combining (\ref{Zt}) and (\ref{NewYPi})
\beqs
|\delta Z^\pi_{t_i}|  &\leq& \left| \mathbb{E} \left[\mathcal{E}_{t_i,T}\partial_x g(X_T^\eps)D_{t_i} X_T^\eps - \mathcal{E}^\pi_{t_{i+1}, t_n}\partial_x g(X^\pi_T) D_{t_i}X^\pi_T \Big/ \mathcal{F}_{t_i} \right] \right| \\
            &=& \left|\mathbb{E} \left[ \mathcal{E}_{t_i,T}\Big(\Big[ \partial_xg(X_T^\eps) - \partial_xg(X^\pi_T)\Big]D_{t_i} X^\pi_T + \partial_xg(X_T^\eps)\Big[D_{t_i} X_T^\eps- D_{t_i} X^\pi_T\Big]\Big)  \right.\right.\\
            &&  \s  + \left.\left.\partial_x g(X^\pi_T)D_{t_i} X^\pi_T \Big[\mathcal{E}_{t_i,T} - \mathcal{E}^\pi_{t_i,T}\Big]  \Big/ \mathcal{F}_{t_i} \right]\right|.
\enqs
From Lemma \ref{LemDXpi}, inequality (\ref{DXXpi}) and Assumption (A2)
\beqs
\mathbb{E} \sup_{0 \leq i \leq n}|\delta Z^\pi_{t_i}|^p  &\leq& \left[\mathbb{E}\left( \sup_{0 \leq i \leq n}D_{t_i}X^\pi_T\right)^\frac{p}{p-1}\right]^{p-1}\Big[\mathbb{E}\Big( \partial_xg(X^\eps_T) - \partial_xg(X^\pi_T)\Big)^p\Big] \\
        &&+ \left[\mathbb{E}\Big( \partial_xg(X^\eps_T)\Big)^\frac{p}{p-1}\right]^{p-1} \left[\mathbb{E}\left( \sup_{0 \leq i \leq n}|D_{t_i} X^\eps_T- D_{t_i} X^\pi_T|^p\right)\right]\\
        &&+ \mathbb{E} \left[ \sup_{0 \leq i \leq n} \mathbb{E}\left[\partial_x g(X^\pi_T)D_{t_i} X^\pi_T (\mathcal{E}_{t_i,T} - \mathcal{E}^\pi_{t_i,T}) \Big/ \mathcal{F}_{t_i} \right]^p\right] \\
        &\leq& C_p\left( \pi^{p/2}+ \mathbb{E}\sup_{0 \leq i \leq n}|I_i|^p \right).
\enqs \\
Using the fact that $|e^x - e^y| \leq (e^x + e^y) |x-y|$, leads to \\
\beqs
I_i   &\leq& C \mathbb{E} \Big\{ \left( D_{t_i} X^\pi_T \mathcal{E}_{t_i,T} + \mathcal{E}^\pi_{t_{i+1}, t_n} \right)  \\
    && \times \left|\int_{t_i}^T \left[f_2(r) -  \frac{1}{2}\int_{E_\eps}[f_3(r)]^2\rho^2(e)m(de) \right]dr \right. + \int_{t_i}^T\int_{E_\eps}f_3(r)\rho(e)\bar{M}(de, dr)  \\
    &&- \sum_{k=i+1}^{n-1}\int_{t_k}^{t_{k+1}}\int_{E_\eps}f_3(r)\rho(e)\bar{M}(de, dr)  \\
    &&  - \left.  \left. \sum_{k=i+1}^{n-1}\int_{t_k}^{t_{k+1}}  \Big[f_2(r) + \frac{1}{2}\int_{E_\eps}[f_3(r)]^2\rho^2(e)m(de) \Big]dr \right| \Big/ \mathcal{F}_{t_i} \right\}\\
\enqs
It follows that
\beqs
I_i    &\leq& C \mathbb{E} \Big\{ D_{t_i} X^\pi_T\left( \mathcal{E}_{t_i,T} + \mathcal{E}^\pi_{t_{i+1}, t_n} \right)  \\
    && \times \left|\int_{t_i}^{t_{i+1}} \left[f_2(r) -  \frac{1}{2}\int_{E_\eps}[f_3(r)]^2\rho^2(e)m(de) \right]dr + \int_{t_i}^{t_{i+1}}\int_{E_\eps}f_3(r)\rho(e)\bar{M}(de, dr) \right| \Big/ \mathcal{F}_{t_i} \Big\}.
\enqs
From Assumption (A.3), we have
\begin{eqnarray*}
|D_{t_i} X^\pi_T \mathcal{E}_{t_i,T}| &\leq& |D_{t_i} X^\pi_T |^r exp\left\{\int_{t_i}^T f_2(u)  du- \frac{1}{2}\int_{t_i}^T\int_{E_\epsilon}f_3(u)\rho^2(e)m(de) \right.du \\
                    && \s \s  \s\s + \left.\int_{t_i}^T\int_{E_\epsilon}f_3(u)\rho(e)\bar{M}(de, du)  \right\} \\
                    &\leq& C \left(\sup_{0 \leq \theta \leq T}|D_{\theta} X^\pi_T |\right)\left( \sup_{0 \leq t \leq T}exp\left\{\int_t^T\int_{E_\epsilon} f_3(t)\rho(e)\bar{M}(de, du)  \right\}\right).
\end{eqnarray*}
Similarly, we have:
\begin{eqnarray*}
|D_{t_i} X^\pi_T |\mathcal{E}^\pi_{t_i,t_n} &\leq& C \left(\sup_{0 \leq \theta \leq T}|D_{\theta} X^\pi_T |\right)\left( \sup_{0 \leq t \leq T}exp\left\{\int_t^T\int_{E_\epsilon} f_3(t)\rho(e)\bar{M}(de, du)  \right\}\right).
\end{eqnarray*}
From other side, for any $r\geq 0 $, we have by Hölder inequality and Proposition \ref{propE}:
\beqs
&&E\left(\sup_{0 \leq t \leq T} exp\left({\int_{t}^T\int_{E_\eps}f_3(r)\rho(e)\bar{M}(de, dr)}\right)\right)^r \\
               && \s \s \leq   \mathbb{E} \left(exp\left\{2r\int_0^T\int_{E_\eps}f_3(r)\rho(e)\bar{M}(de, dr)\right\}\right)^\frac{1}{2} \\
               &&\s \s\s \times \mathbb{E}\left(\sup_{0 \leq t \leq T} exp\left\{-2r\int_0^t\int_{E_\eps}f_3(r)\rho(e)\bar{M}(de, dr)\right\}\right)^\frac{1}{2} \\
               && \s \s \leq   \mathbb{E} \left(exp\left\{r^2\int_0^T\int_{E_\eps}f^3(r)\rho^2(e)m(de) dr\right\}\right)\\
               &&\s \s\s \times \mathbb{E}\left(\sup_{0 \leq t \leq T} exp\left\{-2r\int_0^t\int_{E_\eps}f_3(r)\rho(e)\bar{M}(de, dr)\right\}\right)^\frac{1}{2}\\
               && \s \s <\infty.
\enqs
Thus, for $ p' \in (p, \frac{q}{2})$
\begin{eqnarray*}
\E \sup_{0 \leq i \leq n-1}I^p_i &\leq& C \mathbb{E} \left(\left(\sup_{0 \leq \theta \leq T}|D_{t_i} X^\pi_T |\right)^p \left( \sup_{0 \leq t \leq T}exp\left\{\int_t^T\int_{E_\epsilon} f_3(t)\rho(e)\bar{M}(de, du)  \right\}\right)^p\right. \\
    && \times \left[\sup_{0\leq i \leq n-1}\int_{t_i}^{t_{i+1}} |f_2(r)|dr  +   \frac{1}{2}\sup_{0\leq i \leq n-1} \int_{t_i}^{t_{i+1}} \int_{E_\eps}[f_3(r)]^2\rho^2(e)m(de) dr \right. \\
    && \s +  \left. \sup_{0\leq i \leq n-1}\left|\int_{t_i}^{t_{i+1}}\int_{E_\eps}f_3(r)\rho(e)\bar{M}(de, dr) \right|\right]^p. \\
\enqs
\beqs
    &\leq& C \left[ \E \left(\sup_{0 \leq \theta \leq T}|D_{t_i} X^\pi_T |\right)^\frac{2pp'}{p'-p}\right]^\frac{p'}{2(p'-p)} \\
    && \times \left[ \E \left( \sup_{0 \leq t \leq T}exp\left\{\int_t^T\int_{E_\epsilon} f_3(t)\rho(e)\bar{M}(de, du)  \right\}\right)^\frac{2pp'}{p'-p}\right]^\frac{p'}{2(p'-p)} \\
    && \times \left[ \E\sup_{0\leq i \leq n-1}\left(\int_{t_i}^{t_{i+1}} |f_2(r)|dr \right)^{p'}\right.\\
    && \s\s +   \E\sup_{0\leq i \leq n-1} \left(\int_{t_i}^{t_{i+1}} \int_{E_\eps}[f_3(r)]^2\rho^2(e)m(de) dr \right)^{p'}  \\
    && \s\s +  \left. \E\sup_{0\leq i \leq n-1}\left|\int_{t_i}^{t_{i+1}}\int_{E_\eps}f_3(r)\rho(e)\bar{M}(de, dr) \right|^{p'}\right]^\frac{p}{p'}\\
    &\leq& C \Big[I_1 + I_2 + I_3  \Big]^\frac{p}{p'}.
\end{eqnarray*}
We first estimate $I_3$. \\By Holder inequality for $r>1$, Jensen inequality and Burkholder-Davis-Gundy inequality:
\beqs
 I_3^\frac{p}{p'} &=& \left[\E\sup_{0\leq i \leq n-1}\left|\int_{t_i}^{t_{i+1}}\int_{E_\eps}f_3(r)\rho(e)\bar{M}(de, dr) \right|^{rp'}\right]^\frac{p}{rp'}\\
                &\leq& \mathbb{E} \left(\sum_{0 \leq i \leq n-1}\left|\int_{t_i}^{t_{i+1}}\int_{E_\eps}f_3(r)\rho(e)\bar{M}(de, dr) \right|^{rp'} \right)^\frac{p}{rp'}\\
                &\leq& C\left(\sum_{0 \leq i \leq n-1}\mathbb{E} \left|\int_{t_i}^{t_{i+1}}\int_{E_\eps}[f_3(r)]^2\rho^2(e)m(de) dr\right|^\frac{rp'}{2} \right)^\frac{p}{rp'} \\
                &\leq& C|\pi|^{\frac{p}{2}-\frac{p}{rp'}}.
\enqs
For $\pi$ small enough, we take $r=\frac{2log\frac{1}{|\pi|}} {p'}$, then
\beqs I_3^\frac{p}{p'} &\leq&  C|\pi|^{\frac{p}{2}-\frac{p}{2log\frac{1}{|\pi|}}}. \enqs
And
\beqs
\mathbb{E}\left[I_1 + I_2 \right]^\frac{p}{p'} &\leq& C|\pi|^p.
\enqs
Consequently,
\begin{eqnarray}\label{Zi}
\mathbb{E} \sup_{0 \leq i \leq n}|\delta Z^\pi_{t_i}|^p   &\leq& C|\pi|^{\frac{p}{2}-\frac{p}{2log\frac{1}{|\pi|}}}.
\end{eqnarray}
{\bf Step 2.}
We show that \beqs \E \sup_{0 \leq i \leq n}|\delta \Gamma^\pi_{t_i}|^p  \leq C |\pi|^{p-1},\enqs where $ \delta \Gamma^\pi_{t_i} = \Gamma_{t_i} -\Gamma^\pi_{t_i}$.
In fact,
\beqs
\mathbb{E} \sup_{0 \leq i \leq n}|\delta \Gamma^\pi_{t_i}|^p &=& \mathbb{E} \sup_{0 \leq i \leq n}\left|\int_{E_\eps} \rho(e)\delta U^\pi_{t_i,e} \nu(de)\right|^p \\
            &\leq& C\mathbb{E} \sup_{0 \leq i \leq n}\int_{E_\eps} \rho^p(e)\left|\delta U^\pi_{t_i,e}\right|^p \nu(de) \\
            &\leq& C \mathbb{E} \sup_{0 \leq i \leq n}\left|\delta U^\pi_{t_i,e}\right|^p
\enqs
However, following exactly the same arguments as \emph{Step 1}, we can prove that:
\beqs
\mathbb{E} \sup_{0 \leq i \leq n}\left|\delta U^\pi_{t_i,e}\right|^p \leq C |\pi|^{\frac{p}{2}-\frac{p}{2log\frac{1}{|\pi|}}},
\enqs
and that
\begin{eqnarray}\label{Gammai}
\mathbb{E} \sup_{0 \leq i \leq n}\left|\delta \Gamma_{t_i}\right|^p \leq C|\pi|^{\frac{p}{2}-\frac{p}{2log\frac{1}{|\pi|}}}.
\end{eqnarray}
{\bf Step 3.}: We show that \beqs\mathbb{E} \sup_{0 \leq i \leq n}|\delta Y^\pi_{t_i}|^p  \leq C |\pi|^{p-1}.\enqs
We have
\beqs
Y^\pi_{t_i}&=&\mathbb{E}\left( g(X^\pi_T) + \sum_{k=i+1}^{n-1}f(\Theta^\pi_{t_{k+1}})\Delta t_k/ \mathcal{F}_{t_i}\right),\\
Y^\eps_{t_i}&=&\mathbb{E}\left( g(X_T^\eps) + \sum_{k=i+1}^{n-1}f(\Theta_{t_{k+1}})\Delta t_k/ \mathcal{F}_{t_i}\right).\\
\enqs
Hence by adapted again the argument of \cite{HNS11} to our setting, we have for $ i= n-1, n-2, ...,0.$
\beqs
|\delta Y^\pi_{t_i}| \leq \mathbb{E} \left( \sum_{k=i+1}^{n-1}|f(\Theta^\pi_{t_{k+1}}) - f(\Theta^\eps_{t_{k+1}})|\Delta t_k + |R^\pi_{t_i}| + |\delta g^\pi(X_T^\eps)|\Big/\mathcal{F}_{t_i} \right)
\enqs
where $|\delta g^\pi(X_t^\eps)| = |g(X_T^\eps) - g(X^\pi_T)|$ and
\beqs
|R^\pi_{t_i}| &=& \left|\int_t^T f(\Theta^\eps_r)dr - \sum_{k=i+1}^{n-1}|f(\Theta^\pi_{t_{k+1}})\Delta t_k \right|.
\enqs
For $ j= n-1, n-2, ..., i$
\beqs
|\delta Y^\pi_{t_j}| \leq \mathbb{E} \left( \sum_{k=i+1}^{n-1}|f(\Theta^\pi_{t_{k+1}}) - f(\Theta_{t_{k+1}})|\Delta t_k + \sup_{0 \leq t \leq T}|R^\pi_{t_i}| + |\delta g^\pi(X^\eps_T)|\Big/\mathcal{F}_{t_j} \right)
\enqs
Since we know from \cite{HNS11} that $\mathbb{E}\sup_{0 \leq t \leq T}|R^\pi_{t_i}|^p \leq C|\pi|^\frac{p}{2}$, combining this with standard estimate of $\delta X^\pi$
and Lipschitz property of the generator $f$ we have:
\beqs
&&\mathbb{E}\sup_{0 \leq j \leq n}|\delta Y^\pi_{t_j}|^p \\
&&\leq C\mathbb{E} \left[\left( \sum_{k=i+1}^{n-1}|f(\Theta^\pi_{t_{k+1}}) - f(\Theta_{t_{k+1}})|\Delta t_k \right)^p+ \sup_{0 \leq t \leq T}|R^\pi_{t_i}|^p + |\delta g^\pi(X_T)|^p  \right] \\
&&\leq C\mathbb{E} \left[\left( \sum_{k=i+1}^{n}|\delta X^\pi_{t_k}|\Delta t_k \right)^p+ \left( \sum_{k=i+1}^{n}|\delta Y^\pi_{t_k}|\Delta t_k \right)^p+ \left( \sum_{k=i+1}^{n}|\delta \Gamma^\pi_{t_k}|\Delta t_k \right)^p \right.\\
&& \  \ \ + \left.\sup_{0 \leq t \leq T}|R^\pi_{t_i}|^p + |\delta g^\pi(X_T)|^p\right] \\
&&\leq C\left\{ (T-t_i)^p \mathbb{E} \sup_{i+1 \leq k \leq T} |\delta Y^\pi_{t_k}|^p + \left(|\pi|^{\frac{p}{2}-\frac{p}{2log\frac{1}{|\pi|}}} + |\pi|^\frac{p}{2}\right) \right\}
\enqs
Using similar recursive methods of Theorem 4.2 in \cite{HNS11} we get estimate:
\begin{eqnarray}\label{Y_i}
\mathbb{E}\sup_{0 \leq j \leq n}|\delta Y^\pi_{t_j}|^p\leq C_p|\pi|^{\frac{p}{2}-\frac{p}{2log\frac{1}{|\pi|}}}.
\end{eqnarray}
Finally, combining (\ref{Zi})-(\ref{Gammai}) and (\ref{Y_i}) we close this proof.
\ep\vspace{5mm}\\
{\bf Proof. of Theorem \ref{ThErrorYEpsi}. }\vspace{3mm}\\
For any $i = n-1, n-2, ..., 0.$ Observing that
\beqs
\max_{0 \leq i \leq n} \sup_{t\in[t_i, t_{i+1}]} \mathbb{E}  \left[|Y^\eps_t -Y^\pi_{t_i}|^2\right] &\leq&
                        C\max_{0 \leq i \leq n} \sup_{t\in[t_i, t_{i+1}]} \mathbb{E}  \Big[|Y^\eps_t -Y_t|^2 +|Y_t - Y_{t_i}|^2 +|Y_{t_i}- Y^\pi_{t_i}|^2\Big].
\enqs
Combining  (\ref{errorYEpsl}), (\ref{ErrorYi}) with Corollary 2.7 applied to our setting, we have:
\begin{eqnarray}\label{1}
\max_{0 \leq i \leq n} \sup_{t\in[t_i, t_{i+1}]} \mathbb{E}  \left[|Y^\eps_t -Y^\pi_{t_i}|^2\right]&\leq& C \left(\sigma(\eps)^2 + |\pi|^{1-\frac{1}{log\frac{1}{|\pi|}}}  \right) \label{maxsupY}.
\end{eqnarray}
Combining  (\ref{errorYEpsl}) and (\ref{Zi}) we obtain
\vspace{2mm}
\beq\E\left|\int_0^TV_rdR_r  -\sum_{i=0}^{n-1}Z^\pi_{t_i}\Delta W_{t_i}\right|^2  &\leq& C\left(\E\left|\int_0^TV_rdR_r - \int_0^TZ^\eps_rdW_r\right|^2 \right. \nn \\
            & &   \left.\s+ \E\left|\sum_{i=0}^{n-1}  \int_{t_i}^{t_{i+1}}\Big[Z^\eps_r - Z^\pi_{t_i}\Big]dW_r\right|^2 \right) \nn\\
            &\leq& C\Big( \sigma(\eps)^2 + |\pi|^{1-\frac{1}{log\frac{1}{|\pi|}}}   \Big)\label{2}.
\enq
Arguing as above, we obtain
\beqs
\sum_{i=0}^{n-1}\int_{t_i}^{t_{i+1}}\mathbb{E}|\Gamma_{t} -\Gamma^\pi_{t_i}|^2 dt &\leq& \sum_{i=0}^{n-1}\int_{t_i}^{t_{i+1}}\mathbb{E}\left[|\Gamma_{t} - \Gamma^\eps_t|^2 + |\Gamma^\eps_t - \Gamma^\eps_{t_i}|^2  + |\Gamma^\eps_{t_i} - \Gamma^\pi_{t_i}|^2 \right]dt.
\enqs
From (\ref{errorYEpsl}), (\ref{Gt-Gs}) and (\ref{Gammai}) we have:
\beq
\sum_{i=0}^{n-1}\int_{t_i}^{t_{i+1}}\mathbb{E}|\Gamma_{t} -\Gamma^\pi_{t_i}|^2 dt &\leq& \int_0^T\mathbb{E}\left[|\Gamma_{r} - \Gamma^\eps_r|^2 \right]dr +  C|\pi|^{1-\frac{1}{log\frac{1}{|\pi|}}} \nonumber \\
                    &\leq& C \left(\int_0^T\int_{E_\eps} \E|V_r - U_r^\eps (e)|^2 m(de)dr + \int_0^T\int_{E^\eps} \E|V_r |^2 m(de) dr  \right. \nn\\
                    &&  \left.\s\s +|\pi|^{1-\frac{1}{log\frac{1}{|\pi|}}}\right) \nn \\
                    &\leq& C \left( \sigma(\eps)^2 +  |\pi|^{1-\frac{1}{log\frac{1}{|\pi|}}}   \right) \label{3}.
\enq
Combining (\ref{1}), (\ref{2}) and (\ref{3}) we get the required results.
\ep



\section{Appendix: A priori estimates}
\setcounter{equation}{0} \setcounter{Assumption}{0}
\setcounter{Theorem}{0} \setcounter{Proposition}{0}
\setcounter{Corollary}{0} \setcounter{Lemma}{0}
\setcounter{Definition}{0} \setcounter{Remark}{0}

{\bf Proof.} \textbf{of Lemma \ref{LemDXpi}} \\
By Jensen inequality
\beqs
    \sup_{\theta \leq s \leq t}\|D_\theta X^\pi_t\|^{2p} &\leq & 3^{2p-1} \sup_{\theta \leq s \leq t}\| \int_\theta^s \partial_xb(X^\pi_{\phi^n_r})D_\theta X^\pi_{\phi^n_r}dr\|^{2p} \\
                        &+&   3^{2p-1}\sup_{\theta \leq s \leq t}\| \sigma(\eps)\int_\theta^s \partial_x\beta(X^\pi_{\phi^n_r})D_\theta X^\pi_{\phi^n_r}dW_r\|^{2p} \\
    &+& 3^{2p-1}\sup_{\theta \leq s \leq t}\| \int_\theta^s \int_{E_\eps}\partial_x \beta(X^\pi_{\phi^n_r})D_\theta X^\pi_{\phi^n_r}\bar{M}(de,dr)\|^{2p} \\
                    &+& 3^{2p-1}| \sigma(\eps)\beta(X^\pi_\theta)|^{2p},
\enqs
Taking expectation in both hand side and using Hölder inequality, we obtain:
\beqs
    \mathbb{E}\left[\sup_{\theta \leq s \leq t}\|D_\theta X^\pi_t\|^{2p} \right] &\leq & C_p\left( \mathbb{E}\left[\int_\theta^t \| D_\theta X^\pi_{\phi^n_r}\|^{2p}dr\right] + \mathbb{E}\left[\int_\theta^t \| D_\theta X^\pi_{\phi^n_r}\|^2 dr\right]^p\right.\\
                    &&\s+ \mathbb{E}\left[\int_0^T  \|D_\theta X^\pi_{\phi^n_r}\|^{2p}dr \right] +  \mathbb{E}\Big[\left.| \sigma(\eps)\beta(X^\pi_\theta)|^{2p} \Big]\right)\\
                    &\leq& C\left( \mathbb{E}\left[\int_\theta^t \| D_\theta X^\pi_{\phi^n_r}\|^{2p}dr\right] + B \right)\\
                    &\leq& C\left( \mathbb{E}\left[\int_\theta^t \sup_{0 \leq u \leq r} \| D_\theta X^\pi_{\phi^n_u}\|^{2p}dr\right] + B \right),
\enqs
where $B:=\mathbb{E}\left[\sup_{0\leq s \leq T}\| \beta(X^\pi_s)\|^{2p} \right] $.
Since the constant $C$ doesn't depends on $\theta $ and $n$, we conclude by Gronwall lemma that
$\mathbb{E}\left[\sup_{\theta \leq t \leq T}\|D_\theta X^\pi_t\|^{2p} \right]$ is bounded and therefore $  \sup_{0\leq \theta \leq T}\sup_{n \geq 1 } \mathbb{E}\left[ \sup_{\theta \leq t \leq T}\| D_\theta X^\pi_t\|^{2p} \right] $ is finite . \2 \\
By the same arguments, we prove that $$\sup_{0\leq \theta \leq T}\sup_{n \geq 1 } \mathbb{E}\left[ \sup_{\theta \leq t \leq T}\| D_{\theta,e} X^\pi_t\|^{2p} \right]<  \infty.$$
\ep

\begin{Lemma}
\label{LemAprioriEstimateBSDEs}
Let $\xi \in L^q(\Omega)$,  $f: \Omega\times [0, T] \times \mathbb{R}\times L^2(E, \mathcal{E}, \nu; \mathbb{R}) \rightarrow \mathbb{R}$ be
$ \mathcal{P} \times \mathcal{B} \times \mathcal{B}(L^2(E, \mathcal{E}, \nu; \mathbb{R}))$ measurable, satisfies   $\mathbb{E}\int_0^T|f(t, 0, 0)|^2dt < \infty$ and uniformly Lipschitz w.r.t $(y, z)$, such that, for some constant $ K > 0$ we have:

            $$ |f(t, y_1, u_1) - f(t, y_2, u_2)| \leq K (|y_1- y_2| + \|u_1- u_2\| ),$$
for all $y_1, y_2,  \in \mathbb{R}$ and $u_1, u_2 \in L^2(E, \mathcal{E}, \nu; \mathbb{R})$. \2  \\
Then, there exist a unique triple $(Y, Z, U) \in \mathcal{B}^2$ solution to BSDEs (\ref{FBSDE_Eps}). Moreover, For $q \geq 2$, we have the following a priori estimate:
\begin{eqnarray} \label{Y^q}
\mathbb{E}   \left[ \sup_{0 \leq t \leq T} |Y_t|^q \right]  &+& \mathbb{E} \left(  \int_0^T|Z_t|^2   \right)^{\frac{q}{2}} + \mathbb{E} \left( \int_0^T \int_E |U_t(e)|^q \nu(de)dt \right) \nonumber \\
            & \leq &  C \left(  \mathbb{E} |\xi|^q + \mathbb{E}\left( \int_0^T |f(t, 0, 0)|^2dt \right)^{\frac{q}{2}}  \right).
\end{eqnarray}
\end{Lemma}
{\bf Proof.} Existence and uniqueness of solution of BSDEs with jump are proved in \cite{BBP97} and the estimate (\ref{Y^q}) is a direct consequence of proposition
2.2 in the same reference, with $(f', Q')=(0,0)$.
\ep

%
%
%
%


\end{document}